# Continuous Order Identification of PHWR Models Under Step-back for the Design of Hyper-damped Power Tracking Controller with Enhanced Reactor Safety

Saptarshi Das[a,b], Sumit Mukherjee[b,c], Shantanu Das[d], Indranil Pan[b] and Amitava Gupta[a,b]

*a) School of Nuclear Studies and Applications, Jadavpur University, Salt Lake Campus, LB-8, Sector 3, Kolkata-700098, India.*
*b) Department of Power Engineering, Jadavpur University, Salt Lake Campus, LB-8, Sector 3, Kolkata-700098, India.*
*c) Department of Electrical, Computer and Systems Engineering, Rensselaer Polytechnic Institute, Troy, New York-12180, United States.*
*d) Reactor Control Division, Bhabha Atomic Research Centre, Mumbai-400085, India.*

**Authors' emails:**
saptarshi@pe.jusl.ac.in, s.das@soton.ac.uk (S. Das*)
mukhes3@rpi.edu (S. Mukherjee)
shantanu@magnum.barc.gov.in (Sh. Das)
indranil@pe.jusl.ac.in, indranil.jj@student.iitd.ac.in (I. Pan)
amitg@pe.jusl.ac.in (A. Gupta)

**Abstract:**
In this paper, discrete time higher integer order linear transfer function models have been identified first for a 500 MWe Pressurized Heavy Water Reactor (PHWR) which has highly nonlinear dynamical nature. Linear discrete time models of the nonlinear nuclear reactor have been identified around eight different operating points (power reduction or step-back conditions) with least square estimator (LSE) and its four variants. From the synthetic frequency domain data of these identified discrete time models, fractional order (FO) models with sampled continuous order distribution are identified for the nuclear reactor. This enables design of continuous order Proportional-Integral-Derivative (PID) like compensators in the complex *w*-plane for global power tracking at a wide range of operating conditions. Modeling of the PHWR is attempted with various levels of discrete commensurate-orders and the achievable accuracies are also elucidated along with the hidden issues, regarding modeling and controller design. Credible simulation studies are presented to show the effectiveness of the proposed reactor modeling and power level controller design. The controller pushes the reactor poles in higher Riemann sheets and thus makes the closed loop system hyper-damped which ensures safer reactor operation at varying dc-gain while making the power tracking temporal response slightly sluggish; but ensuring greater safety margin.

## 1. Introduction

In recent years, fractional order systems, governed by fractional order differential equations have got increased interest in scientific community for the modeling of physical systems [1] with greater accuracy. System identification which is a well established tool in control engineering to build models for unknown or poorly understood dynamical systems has already been extended using fractional calculus to get a better description of the physical system [2]. Application of fractional calculus has also been done in few nuclear engineering problems as a better description of the neutron diffusion equation in spatial [3] and temporal [4] domain, point reactor kinetics [5], neutron transport equation [6]-[7], compact modeling of nuclear reactor [8] and robust controller design for reactor power regulation [9] etc. Fractional order modeling of physical systems firstly requires the knowledge of the number of FO elements, present in the model i.e. the number of terms in the numerator and denominator of the FO transfer function model which optimally describe the dynamical behavior of the system. Then, the fractional orders of the model along with the associated coefficients [10] are estimated from an experimental data-set. Fractional order models whose numerator or denominator orders can be described by a decreasing power of the Laplace variable ($s$) with a simple FO least common divisor, are known as commensurate order models [11]. Also, fractional order models whose orders do not have a least common divisor (having irrational orders or recurring decimal numbers or their truncated version) are known as incommensurate order models. Model reduction of higher order processes in a flexible order template may lead to such incommensurate fractional order models [12]. In fact, the incommensurate order FO models can only be visualized as a commensurate order model with very small commensurate order [11].

The notion of frequency domain modeling of dynamical systems with discrete integer order elements has been extended by Valerio and Sa da Costa [13] to discrete fractional order models. Hartley and Lorenzo [10] first proposed that for a process model, the orders of differentiation do not necessarily have to be discrete in nature. Rather a continuous distribution of orders can be thought of, among which only a limited number of orders are significantly large. These orders can be termed as the dominant fractional orders and can also be represented by the flexible order process model reduction as proposed in [12]. Fractional order model building with time domain system identification techniques are studied in Malti *et al.* [2] with the consideration of noisy measurement. A generalized algorithm for FO system identification with measured frequency domain or time domain data has been proposed by Valerio and Sa da Costa [14]. It has been suggested in [14] that the easiest way to identify FO models is to build higher integer order models with the available system identification techniques in time domain and then generating synthetic frequency domain data out of that model to build compact fractional order models. But in most cases, the experimental data is available in time domain. However for fractional order system identification, most of the robust estimators are developed for frequency domain data. Thus, it was necessary to transform time domain information of the dynamical system into an equivalent frequency domain data. In order to do so, the concept of Valerio and Sa da Costa [14] of higher order time domain system identification and their compact representation using fractional order models have been applied in this paper and the concept is extended wherever needed. Traditionally frequency domain system identification is done by using Levy's method of complex curve fitting, which is a least squares based method that doesn't work equally well at all frequencies. Valerio and Sa da Costa [13] extended this method for commensurate fractional order transfer functions and also improvised the scheme by introducing weights to the basic algorithm and removed the frequency dependence of the method [15]. Vinagre's weights on the Levy's method further enhances the identification methodology but do not always lead to better results as discussed in the literatures [13], [15] with the corresponding software implementation in [16]. As an alternative to data driven system identification and modeling for controller development, the classical ODE/PDE based first principle modeling may be adopted which may also be modeled using FO dynamics e.g. investigation like finite differencing for fractional point kinetics [17], fractional point-kinetics in reactor start-up [18]-[19], time fractional Telegrapher's equation for neutron motion [20]-[21], stochastic point kinetic equation [22], fractional point kinetics for reactor with slab geometry [23] etc. But in all of these cases, accurate knowledge of all parameters of the governing physical equations is essential, which is often impractical for many large physical systems. For the controller design techniques for uncertain FO systems, identification of the structure of uncertainty e.g. additive, multiplicative or interval type etc. is even more difficult in most cases.

The continuous order system identification, proposed by Hartley and Lorenzo [10] is a completely new philosophy of data based system modeling where a continuum in the system's order is considered. In the pioneering work [10], analytical expressions for system's transfer function representation (as transcendental functions of Laplace variable "$s$") have been given for various idealized order-distributions e.g. uniform, Gaussian, triangular, impulsive, truncated ramp type etc. In fact, for continuous order identification of any practical system, the order distributions may not follow these ideal shapes for which closed form analytical expressions exists to represent its transfer function. Discrete/sampled commensurate fractional order system identification and its extension to all pole continuous order system identification has been extensively studied in [10] and its software implementation can be found in the Matlab based toolbox Ninteger [16]. This concept has been extended for frequency domain continuous order system identification with pole-zero models by Nazarian and Haeri [24] with the identifiability conditions given in [25]. In this paper, the Levy's frequency domain fractional order system identification technique and its improved version with Vinagre's frequency weights [13], [15] have been used with a practical test data and few interesting and new results are also reported. The preset approach considers gradual reduction in the commensurate order of the fractional order model to be fitted with the data while continuously observing its accuracy. In a theoretical sense, when the commensurate fractional order of a model tends towards zero or a very small value, the model can be considered as a continuous order model [10]. We have found that with a finite number of data points, arbitrary reduction in the commensurate order does not always produce a better quality of model, in terms of the modeling error. Rather, for very small commensurate order the number of unknown variables (coefficients of the numerator and denominator) becomes very large and the accuracy becomes poor, due to significant computational errors with large system matrices. For this reason, it is very important to find out an accurate choice of the commensurate order which explains the data correctly, on the other hand intermediate matrices does not become ill-conditioned.

The notion of PID controllers which are widely used in process control has been first extended by Podlubny [26] with the fractional order PID or $PI^{\lambda}D^{\mu}$ controller which has two extra degrees of freedom over the three-term PID controller viz. the integro-differential orders. The fractional order $PI^{\lambda}D^{\mu}$ controller has five independent parameters to tune and takes the following form:

$$C(s)=K_p+\frac{K_i}{s^{\lambda}}+K_d s^{\mu}=\frac{K_d s^{\lambda+\mu}+K_p s^{\lambda}+K_i}{s^{\lambda}} \tag{1}$$

Here, $C(s)$ represents the controller with '$s$'being the Laplace variable or complex frequency. Gains $\{K_p, K_i, K_d\}$ control the mixing of proportional, integral and derivative actions. Integro-differential orders $\{\lambda, \mu\}$ give extra flexibility in balancing the effect of poles and zeros using the concepts of fractional calculus. The concept of three and five term controllers like PID and $PI^{\lambda}D^{\mu}$ respectively, was extended to the generalized continuous order PI/PID controllers by Hartley and Lorenzo [27] which have a continuous distribution of zeros instead of two zeros of a PID controller. The generalized continuous order PI/PID controller takes the form (2) and is expected to give better control performance if it can be tuned properly. Now, generalizing the controller gains $\{K_p, K_i, K_d\}$ in (1) as $\{K_0, K_1, \cdots, K_N\}$ and considering integer order pole with only fractional order zeros we get:

$$C(s)=\frac{K_0 s^{Nq}+K_1 s^{(N-1)q}+\cdots+K_{N-1}s^{q}+K_N}{s}=\frac{\sum_{n=0}^{N}K_n s^{(N-n)q}}{s} \tag{2}$$

In (2), $q$ is the commensurate order of the continuous order PI/PID controller with $q<1, q\in\mathbb{R}_+$ such that $Nq=1$ for PI controller and $Nq=2$ for PID controller respectively. The concept has also been extended in [27] for designing generalized continuous order dynamic compensator for controlling continuous order systems. These controllers have more design flexibility and degrees of freedom as more closed loop poles and zeros can be placed at desired locations by proper selection of its gains unlike placing only two closed loop poles using PID type controllers [28]. Therefore the generalized continuous order compensator takes the form (3) as a further improvement of the scheme in (2). Also, for very small commensurate order ($q\to 0$) the numerator and denominator of (3) can be represented by definite integrals denoting the continuous order distribution for the numerator and denominator of the compensator.

$$C(s)=\frac{\sum_{n=0}^{N}K_n^{num}s^{(N-n)q}}{\sum_{n=0}^{N}K_n^{den}s^{(N-n)q}}\simeq\frac{\int_0^{Nq}K^{num}(q)s^{q}dq}{\int_0^{Nq}K^{den}(q)s^{q}dq} \tag{3}$$

The only problem with the compensator structure (3) is that it lacks the desired set-point tracking capability of PI/PID type controllers due to the absence of an in-build integrator unlike structure (2). Hence, (3) is a generalization of the FO lead-lag compensator introduced in [29]. Therefore, in the present study, the controller design has only been restricted with the structure given in (2). Also, in controller structure (2) we have considered an integer order integrator rather than using a fractional integrator as in (1), because of the fact that former makes the control system work much faster than with latter. Also, in [27] it has been suggested that the identification and controller design methodology can be improved by replacing the summation in the numerator and denominator of the continuous order model and compensator respectively with definite integrals as in (3), thereby considering all possible real orders that may be present along with their corresponding coefficients. However, with this particular method, the resulting controller will be difficult to implement in real hardware due to the constraints involved in realizing the huge number of fractional order operators [30]-[31]. The present paper firstly applies the concept of continuous order identification for a nuclear reactor under step back condition at different operating points and then designs a robust continuous order PID like controller (2) that works at all operating points despite the gravely nonlinear nature of the plant.

Earlier investigations regarding the modeling of operating PHWR under step back [8], [9], show that the dc gain of the nonlinear nuclear reactor gets changed with shift in operating point (initial power and level of control rod drop). PID type controller with fractional order enhancements like FO phase shaper [9] and $PI^{\lambda}D^{\mu}$ controller [8] have been applied to ensure robust operation of the reactor in wide range of operating points. The present paper further enhances the concepts in [8], [9] in the light of continuous order system identification and controller design. It is well known that with the help of classical PID type controllers the dominant closed loop poles of a process can be modified in the complex $s$-plane. For integer order system and controllers, the whole $s$-plane is termed as the primary Riemann sheet. Hartley and Lorenzo [27] have shown that for fractional order systems, the controller design task gets mapped in secondary or tertiary Riemann sheets. The significance of the presence of poles in the higher Riemann sheets can be described as weak non-dominant dynamical behavior of the system. The concept of fractional order systems and control enables the design of pole placement like tuning of process controllers using the possibility of their existence in higher Riemann sheets. This has been found to have extreme importance to doubly ensure safer operation of nuclear reactors. It is well known that the stability of FO systems are more, even in perturbed condition, if all of its poles lie in higher Riemann sheet (hyper-damped or ultra-damped poles). Therefore, the conventional pole placement controller design in $s$-plane can be improved to push all closed loop poles in higher Riemann sheet to achieve higher stability margin. Even in classical integer order controller design, over-damped closed loop poles may exhibit oscillatory response if the process gain is increased heavily due to nonlinearity or any possible mishandling by the operator or under faulty condition. In nuclear reactor power level control such oscillations are strictly prohibited since at low power, the reactor might get poisoned out and the mechanical elements might experience thermal shocks in the presence of oscillating power dynamics. However the proposed design approach has an inherent capability to nicely handle all of these issues. For any fractional order system if the controller is designed with such an objective that all of its closed loop poles lie in the higher Riemann sheet, then to reach instability the branches of the root locus must cross all the secondary and tertiary Riemann sheets and finally the stable region of the primary Riemann sheet. This makes the closed loop system hyper-damped or ultra-damped which can be viewed like extracting much more stability margin than the usual notion of stability for linear integer order systems. But this comes at the cost of sluggish response of the plant although the plant is still able to track the desired set-point if the fractional/continuous order controller contains an integer/fractional order integrator in its structure. This new concept of controller design of pushing all the closed loop poles in the hyper-damped or ultra-damped region doubly ensures the issues like safer reactor operation at varying dc-gain due to nonlinearity though it compromises a bit on the time response. Although higher Riemann sheet poles cause slow time response, such a sacrifice in power level tracking performance is worth to ensure greater reliability and safety features for the control of safety-critical systems like nuclear reactors.

Rest of the paper is organized as follows. Section 2 describes the higher integer order discrete time transfer function modeling from test data of a PHWR under step-back condition. Continuous order modeling approach of the nuclear reactor is described in section 3. Section 4 describes an optimization based pole placement like tuning of continuous order controller to ensure dead-beat power level tracking at wide range of operating points. The paper ends with the conclusion in section 5, followed by the references.

## 2. Discrete time system identification of a nuclear reactor under different step-back conditions

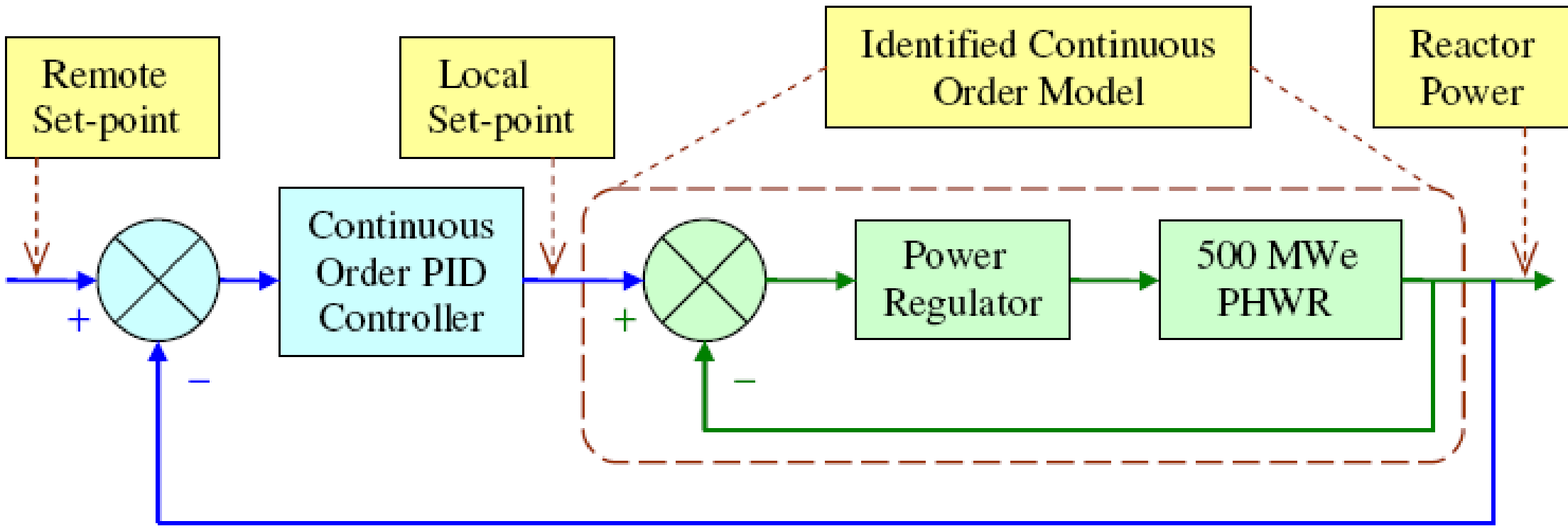

Fig. 1. Identified system comprising of the PHWR along with its power regulator and the proposed modifications in the reactor control scheme.

The major nonlinearity is introduced in the dynamics of a nuclear reactor due to the cross product of state (neutron density) and input (reactivity) in the point kinetic equation. The coupling between neutronics and thermal-hydraulics are almost linear as reported in the investigation by Das *et al.* [32]. The nonlinearity becomes predominant when shift in operating power i.e. initial value of the state and reactivity or equivalent control rod worth i.e. control inputs are varied. This motivates us to study the dynamical behavior of the reactor under four different operating power levels and two different levels of added negativity reactivity. The effect of the nonlinearity is evoked when the reactor power level or added reactivity levels are changed.

In this section, dynamical models are identified for a 500MWe PHWR under various step-back conditions from test-data as studied in [8]-[9]. For this purpose, the reactor needs to be modeled using the dynamics of power variation during a step-back with the change in control rod position as the input and the global reactor power as the controlled variable. The control scheme for the reactor is shown in Fig. 1. Here, the PHWR with its power regulator in closed loop with thermal feedback is taken as the system and is controlled by a master controller as a continuous order PID like compensator with the master controller output acting as the local set point (Fig. 1). In the present modeling scenario, a nonlinear point kinetic equation system has been assumed to describe the nuclear reactor dynamics. Hence the nonlinearity of the system will make the global reactor power transients to differ with each other in a large extent with shift in operating condition i.e. the initial reactor power at which the reactor was operating in steady-state and the level of step-back or negative reactivity addition. Thus the nonlinear dynamical behavior of the nuclear reactor with the power regulator in loop is identified as stable transfer function models around different operating point, though the open loop nuclear reactor model without the power regulator is marginally stable [33]. From these stable higher order discrete time models, continuous order reactor models are developed which enables design of a single continuous order PID like controller that ensures dead-beat power tracking at several operating points.

It is well known that system identification refers to mathematical modeling of dynamical systems where the physics of the process is highly complicated and the system's governing laws are not well understood. It is basically finding an approximate model from an input-output experimental data by an iterative technique, where the modeling requires less insight of the actual system physics. There are several classical identification methods e.g. time response based, frequency response based methods etc. In the present work, a time response based system identification approach is adopted, to find out the transfer function between power developed by a nuclear reactor and the level of control rod drop. For this purpose, the basic least square estimation based system identification techniques and other variants of LSE are briefly introduced next.

### *2.1. Identifying higher order linear models using least square estimator*

Here, the generalized identification technique using recursive-least square algorithm from a measured time domain data [34] is briefly discussed. Let us assume that at time event $t$, the input and output of an unknown system is $u(t)$ and $y(t)$ respectively. Then the system can be described by the following linear difference equation

$$y(t)+a_1 y(t-1)+\cdots+a_n y(t-n)=b_1 u(t-1)+\cdots+b_m u(t-m) \tag{4}$$

The above equation can be re-written in the following form if the values of input and output data at each time step are known

$$y(t)=-a_1 y(t-1)-\cdots-a_n y(t-n)+b_1 u(t-1)+\cdots+b_m u(t-m) \tag{5}$$

The calculated value of the output is thus

$$\hat{y}(t)=\varphi^T(t)\cdot\theta \tag{6}$$

where, system parameters

$$\theta=\left[a_1\cdots a_n\ b_1\cdots b_m\right]^T \tag{7}$$

and measured input-output data

$$\varphi(t)=\left[-y(t-1)\cdots-y(t-n)\ u(t-1)\cdots u(t-m)\right]^T \tag{8}$$

Now, from the input-output data ($Z^N$) over a time interval $\left(1\le t\le N\right)$ the coefficient vector $\theta$ can be calculated satisfying the condition

$$\hat{\theta}=\min_{\theta} V_N(\theta,Z^N) \tag{9}$$

where, $Z^N=\left\{u(1),y(1),\cdots,u(N),y(N)\right\}$ and

$$V_N(\theta, Z^N) = \frac{1}{N}\sum_{t=1}^{N}\left(y(t) - \hat{y}(t|_\theta)\right)^2 = \frac{1}{N}\sum_{t=1}^{N}\left(y(t) - \varphi^T(t)\cdot\theta\right)^2$$

To find out the minimum value in (9), the derivative of $V_N$ with respect to $\theta$ needs to be set to zero.
i.e.

$$\frac{d}{d\theta}V_N(\theta, Z^N) = \frac{2}{N}\sum_{t=1}^{N}\varphi(t)\left(y(t) - \varphi^T(t)\cdot\theta\right) = 0$$

$$\Rightarrow \quad \frac{1}{N}\sum_{t=1}^{N}\varphi(t)y(t) = \frac{1}{N}\sum_{t=1}^{N}\varphi(t)\varphi^T(t)\theta \qquad (10)$$

$$\Rightarrow \quad \hat{\theta}_N^{LS} = \left[\frac{1}{N}\sum_{t=1}^{N}\varphi(t)\varphi^T(t)\right]^{-1}\frac{1}{N}\sum_{t=1}^{N}\varphi(t)y(t)$$

Since in this method, the sum of the squared residuals or errors is minimized, it is known as the least square algorithm for system identification. Also with the known value of input and output data at each instant i.e. $\varphi(t)$ vector, using relation (10) the least square estimate of the coefficients of discrete transfer function model i.e. $\hat{\theta}^{LS}$ can be obtained.

### *2.2. Basic least square estimator and its variants*

The minimization of the identification error depends largely on the structure of the estimator. The choice of a suitable structure for the noise model as well as the system model plays a very important role in minimizing the modeling error. This sub-section briefly describes few variants of basic LSE and their roles in system identification and the choice of a proper estimator structure [35].

Let us consider, a generalized linear model structure of the form

$$y(t) = G(q^{-1}, \theta)u(t) + H(q^{-1}, \theta)e(t) \qquad (11)$$

where, $u(t)$ and $y(t)$ are the input and output of the system respectively, $e(t)$ is the zero-mean white noise, $\theta$ is the parameter vector to be estimated, $G(q^{-1}, \theta)$ is the transfer function of the deterministic part of the system and $H(q^{-1}, \theta)$ is the transfer function of the stochastic part of the system. Here $q^{-1}$ denotes the backward shift operator. Equation (11) can further be rewritten as (12) which is known as the equation error type linear LSE.

$$A(q^{-1})y(t) = \frac{B(q^{-1})}{F(q^{-1})}u(t) + \frac{C(q^{-1})}{D(q^{-1})}e(t) \qquad (12)$$

where, $\{B, F, C \,\&\, D\}$ are polynomial in $q^{-1}$ and represent the numerator and denominator of the system model and noise model respectively and $\{A\}$ represents the polynomial containing common set of poles for both of the system and noise model. The block diagram representation of the generalized least-square estimator is shown in Fig. 2.

The generalized LSE structure (12) can be further customized by considering only fewer elements among $\{B, F, C, D \,\&\, A\}$ at once while choosing different estimators for system identification which are detailed in the following subsections. For example a Finite Impulse Response (FIR) form for the model can be obtained by considering the polynomial $\{B\}$ only etc. The next subsections briefly describe four classes of estimators as special cases of the generalized equation error type linear LSE described by (12).

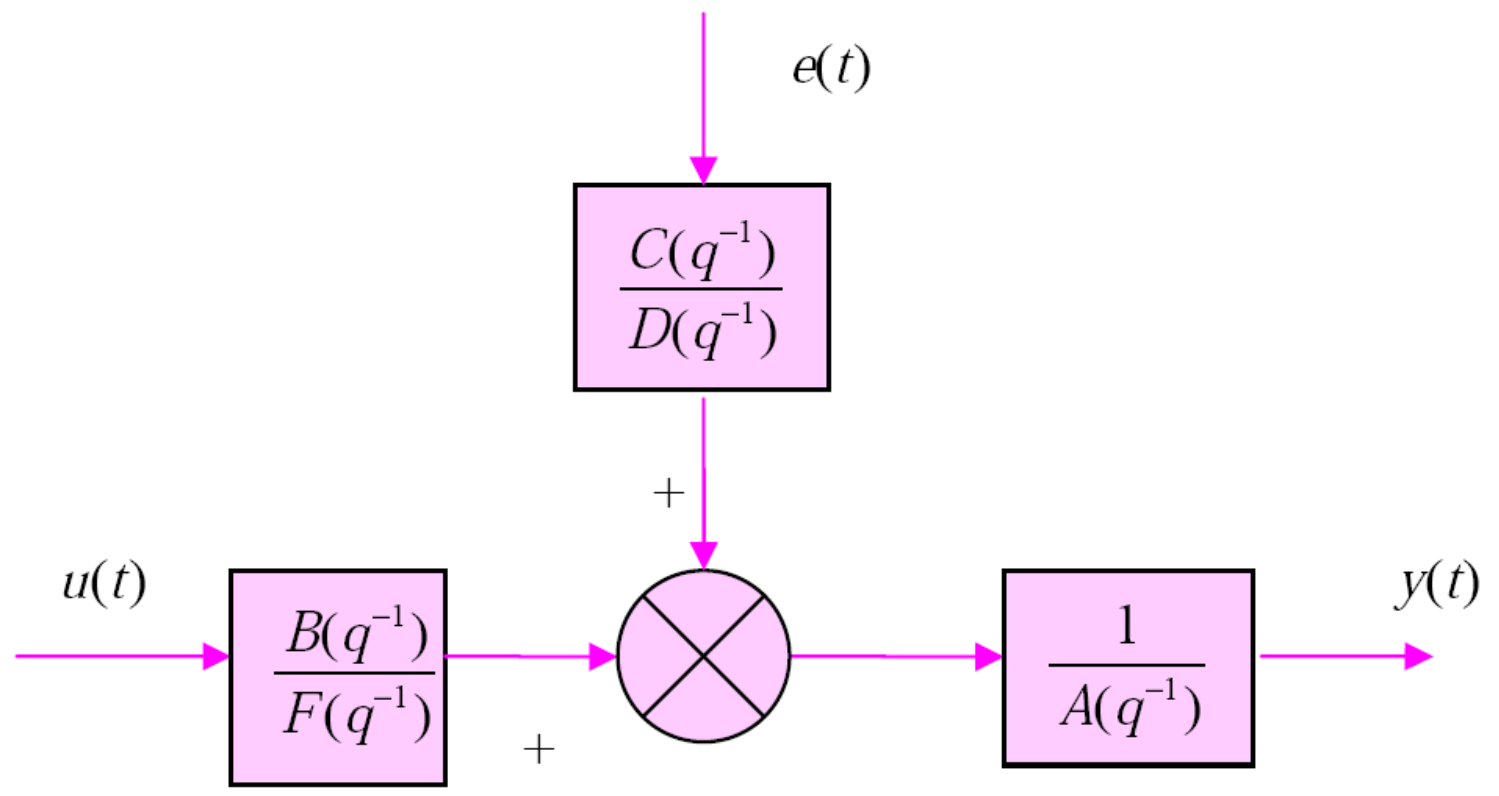

Fig. 2. Block diagram representation for generalized linear least-square estimator.

*2.2.1. AutoRegressive eXogenous (ARX) estimator*

The basic structure of an ARX estimator is governed by (13).

$$A(q^{-1})y(t) = B(q^{-1})u(t) + e(t) \tag{13}$$

This structure doesn't allow modeling of the noise and the system dynamics independently (Fig. 3). The main disadvantage of this structure is that the deterministic (system) dynamics and the stochastic (noise) dynamics are both estimated with same set of poles which may be unrealistic for many practical applications.

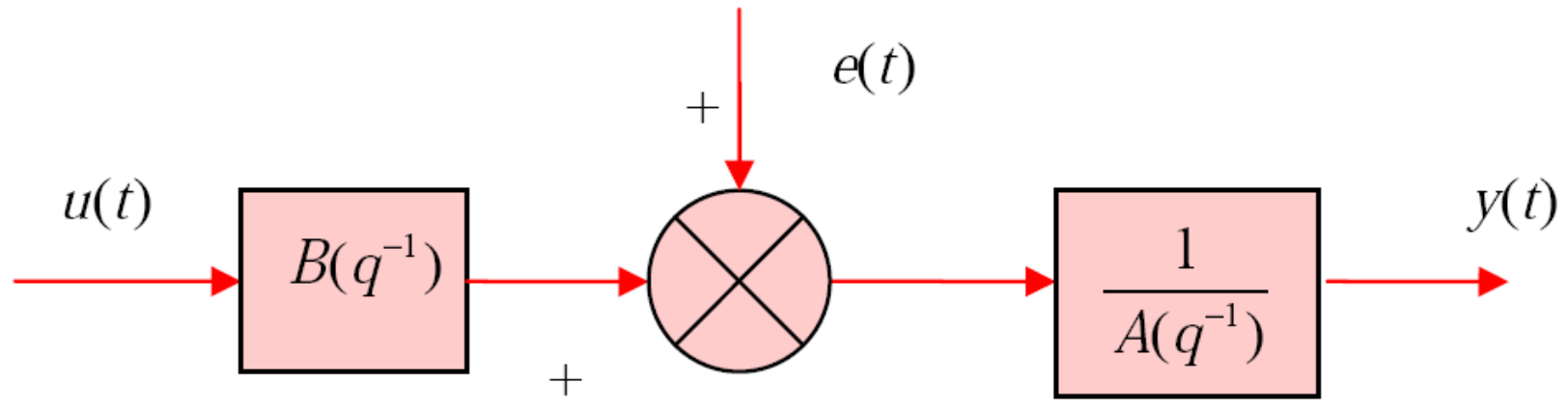

Fig. 3. Block diagram representation for ARX model structure.

*2.2.2. AutoRegressive Moving Average eXogenous (ARMAX) estimator*

Basic structure of an ARMAX estimator is given by (14).

$$A(q^{-1})y(t) = B(q^{-1})u(t) + C(q^{-1})e(t) \tag{14}$$

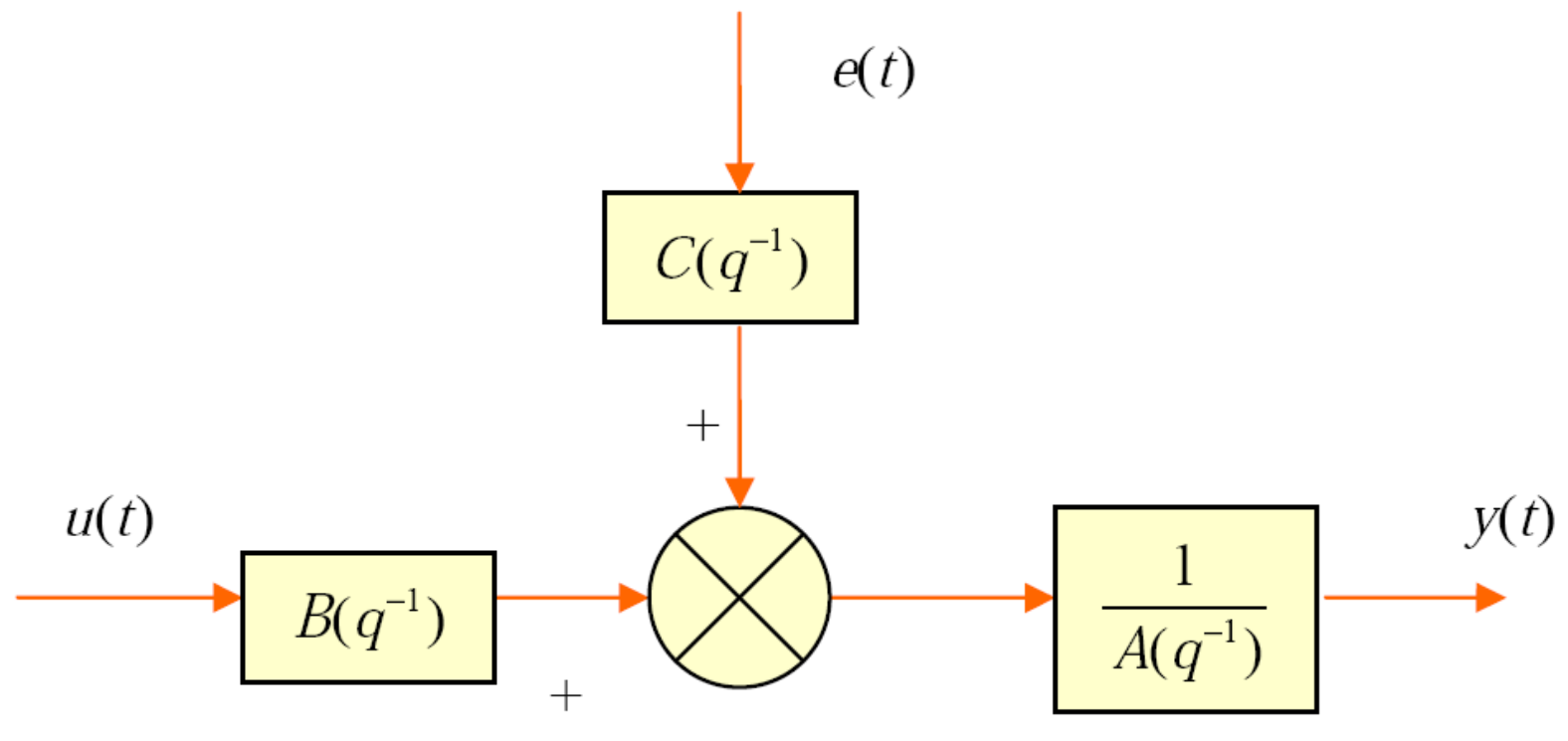

Fig. 4. Block diagram representation for ARMAX model structure.

The basic disadvantage of ARX structure is that it has inadequate freedom to describe the exogenous noise dynamics which could be modeled with better flexibility by introducing a moving average to the white noise. Thus, the ARMAX structure gives better flexibility over ARX structure to model the measurement noise along with the system. ARMAX structure estimates different set of zeros but same set of poles for the system and the noise model (Fig. 4). This structure is especially suitable when the stochastic dynamics are dominating in nature and the noise enters early into the process e.g. load disturbances.

*2.2.3. Box-Jenkins (BJ) estimator*

Basic structure of the BJ estimator is governed by the following relation

$$y(t) = \frac{B(q^{-1})}{F(q^{-1})}u(t) + \frac{C(q^{-1})}{D(q^{-1})}e(t) \tag{15}$$

BJ structure allows estimation of different set of poles and zeros for the system and noise model (Fig. 5). This model structure is especially suitable when disturbances enters into the model at later stage e.g. measurement noise.

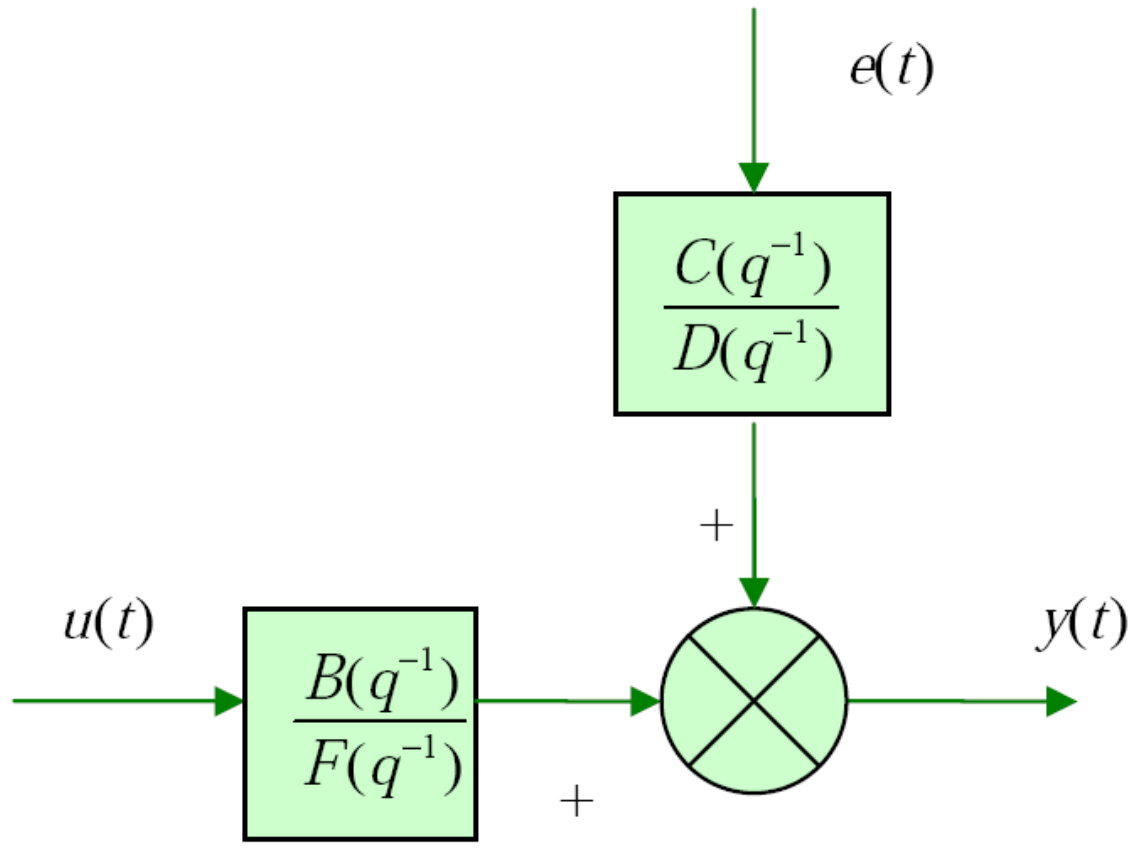

Fig. 5. Block diagram representation for Box-Jenkins model structure.

*2.2.4. Output-Error (OE) estimator*

An OE estimator has the following structure

$$y(t) = \frac{B(q^{-1})}{F(q^{-1})} u(t) + e(t) \quad (16)$$

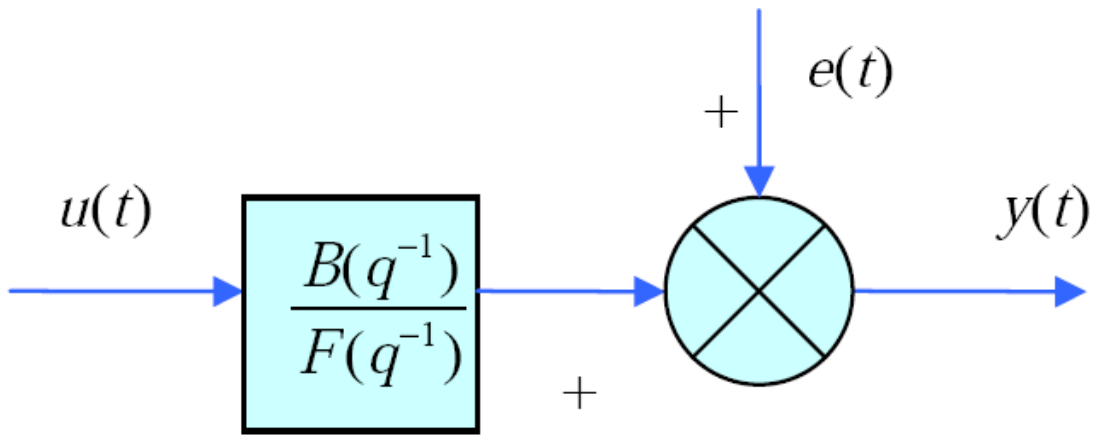

Fig. 6. Block diagram representation for Output-Error model structure.

The OE structure estimates poles and zeros for the system model only. It does not estimate the noise model. This structure is suitable when modeling of the system dynamics is of prior concern and not the noise-model or the measurement noise is negligible.

### *2.3. Time domain identification results and model validation*

The section presents system identification of a PHWR along with its regulating system (Fig. 1) using the above mentioned variants of LSE. For identification, the reactor is visualized as a system with control rod position (fraction of total drop) as input and the global power (in percentage of maximum power produced) as output. The identification is based on data obtained from operating Indian PHWRs provided by Nuclear Power Corporation of India Ltd. (NPCIL) as also studied in [8]-[9]. The data at different step-back levels is provided for 14 seconds with 0.1 second of sampling time. Graphical representation of the data is shown in Fig. 7 for 30% and 50% rod drop cases with different initial powers i.e. 100%, 90%, 80% and 70%. With the data in Fig. 7 (a and b), stable discrete time higher integer order transfer function models are built using the four class of estimators i.e. ARX, ARMAX, BJ and OE as introduced in previous section.

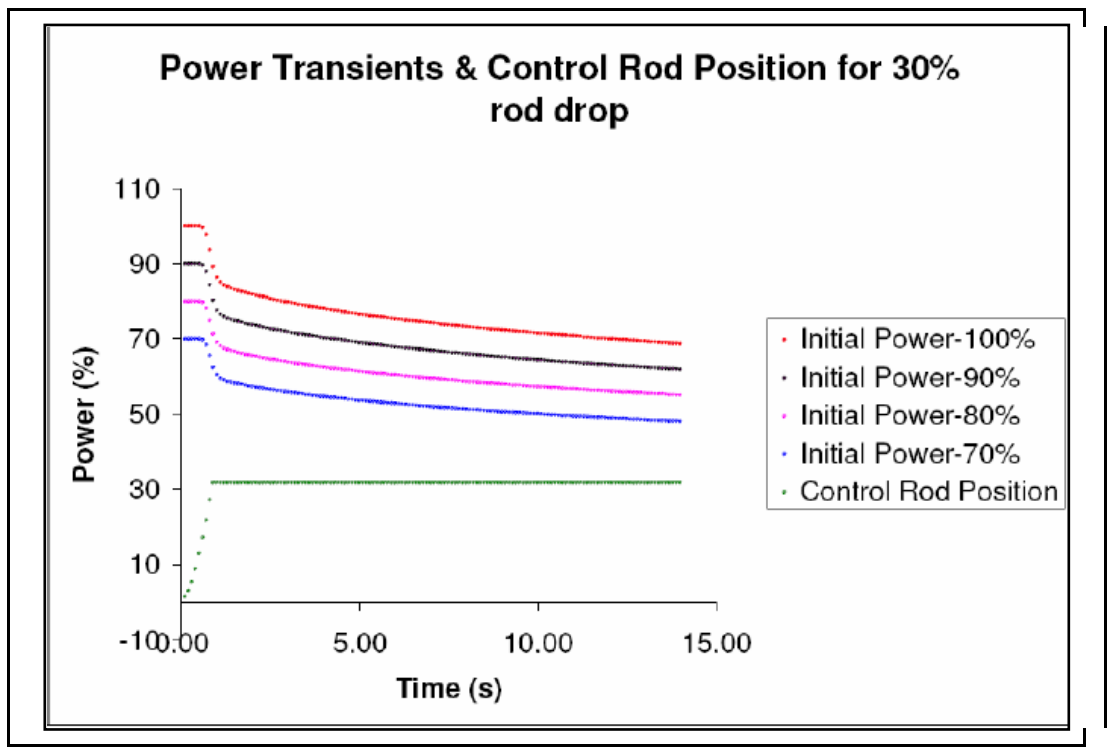

(a)

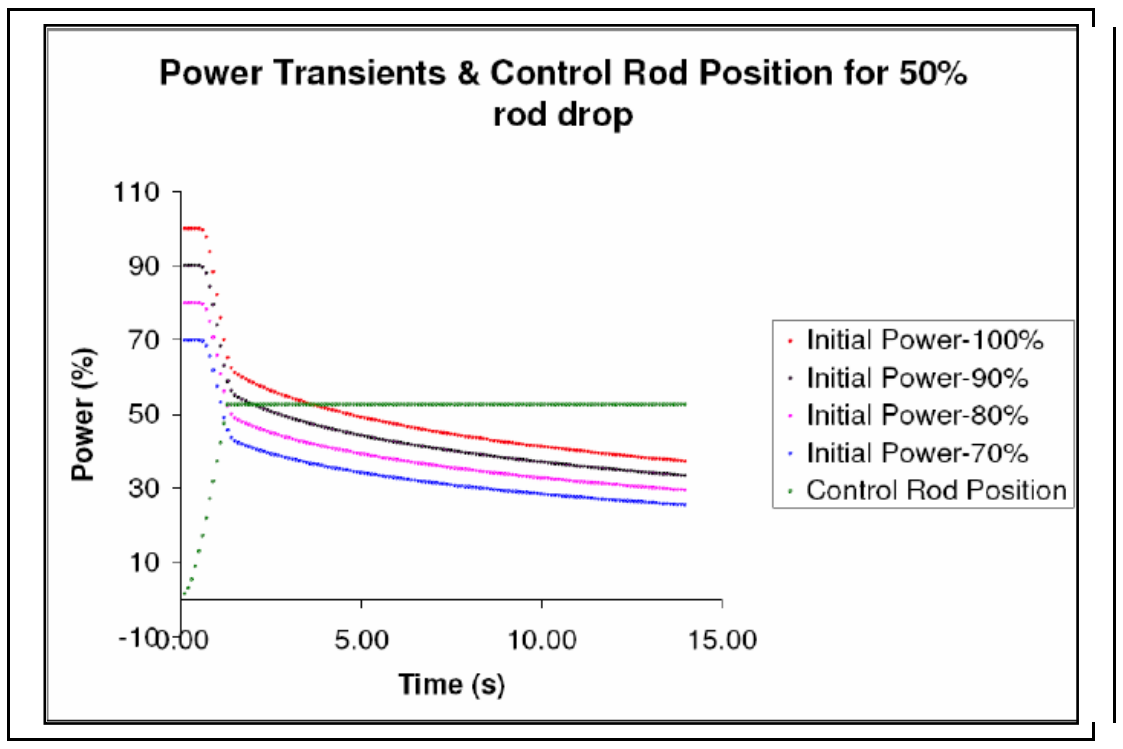

(b)

Fig. 7. Reactor power transient and control rod drop data used for system identification.

Also, it is an essential criterion in system identification that the models should be built in such a manner that it can maximally extract all the information hidden in the data. This capability of a model is judged with the help of few statistical performance criteria like Akaike's Information Criteria (AIC), Final prediction Error (FPE), percentage fit etc. [36]. We have compared the accuracies of the identified discrete-time models with respect to their AIC values which is a common practice in data based modeling of processes and has certain advantages over the other performance criteria [37]. The Akaike's Information Criterion (AIC) [37] is defined as:

$$AIC = \log V + \frac{2d}{N} \tag{17}$$

where $V$ is the loss function, $d$ is the number of estimated parameters, and $N$ is the number of values in the estimation data set.
The loss function $V$ is defined by

$$V = \det\left(\frac{1}{N}\sum_{1}^{N}\varepsilon\left(t,\theta_N\right)\left(\varepsilon\left(t,\theta_N\right)\right)^T\right) \tag{18}$$

where, $\theta_N$ represents the estimated parameters and $\varepsilon$ is the estimation error.

Table 1 reports the best found AIC values of the four classes of estimators with increasing order of estimated models. According to [37], the better model should have a lower AIC value. Also, a trade-off has been made between the significant improvement in the AIC value and the complexity of identified models due to unnecessary increase in system order. It is also observed that even for modeling a nonlinear system around a specific operating point, an increase in system's order does not always result in good modeling performance.

The focus of the paper is to build fractional order models with time domain data of the reactor operation. But firstly we have adopted the approach of identifying higher integer order models. This is because of the fact that regarding time domain data based fractional order model building, the pioneering works like Valerio and Sa da Costa [14] have suggested that fractional order models should be built using the frequency domain information of identified higher integer order discrete time transfer function. In order to do so, the most accurate discrete time models corresponding to the Box-Jenkins estimator in Table 1 at various operating conditions of the reactor are reported in equations (19)-(26). In the identified models ($G$) the superscripts denote the initial reactor power at which step-back is initiated and the subscript denotes the level of control rod drop.

Table 1
Choice of suitable identifier based on minimum modeling error (AIC values)

| **Rod drop level** | **Initial Power** | **System identification algorithm** | | | |
|---|---|---|---|---|---|
| | | **ARX** | **ARMAX** | **Box-Jenkins** | **Output Error** |
| 30% | 100% | -5.5951 | -6.8921 | -6.9892 | -6.986 |
| | 90% | -5.9025 | -7.2617 | -7.2878 | -7.2841 |
| | 80% | -6.0278 | -7.3159 | -7.3438 | -7.2401 |
| | 70% | -6.2675 | -7.3956 | -7.4957 | -7.3775 |
| 50% | 100% | -4.7074 | -6.6619 | -7.0285 | -5.5849 |
| | 90% | -5.1039 | -7.1781 | -7.2645 | -6.6075 |
| | 80% | -5.3236 | -7.271 | -7.3227 | -7.2683 |
| | 70% | -5.4029 | -6.4896 | -6.5131 | -6.3896 |

$$G_{30}^{100}\left(z\right)=\frac{-33.7z^{4}+48.94z^{3}-8.075z^{2}-1.686z-0.7376}{z^{5}-1.485z^{4}+0.5105z^{3}} \quad (19)$$

$$G_{30}^{90}\left(z\right)=\frac{7.773z^{6}}{z^{7}-1.994z^{6}+2.522z^{5}-2.848z^{4}+2.468z^{3}-1.682z^{2}+0.7502z-0.171} \quad (20)$$

$$G_{30}^{80}\left(z\right)=\frac{-18.59z^{6}+29.22z^{5}-11.9z^{4}+16.78z^{3}-7.937z^{2}+3.413z-0.2431}{z^{7}-0.9305z^{6}} \quad (21)$$

$$G_{30}^{70}\left(z\right)=\frac{-30.44z^{6}+48.92z^{5}-14.66z^{4}-3.01z^{3}+7.914z^{2}-4.594z+1.306}{z^{7}-1.409z^{6}+0.5661z^{5}-0.1161z^{4}} \quad (22)$$

$$G_{50}^{100}\left(z\right)=\frac{-0.9878z^{6}}{z^{7}-1.768z^{6}+0.8855z^{5}+0.2743z^{4}-1.02z^{3}+1.157z^{2}-0.6595z+0.1493} \quad (23)$$

$$G_{50}^{90}\left(z\right)=\frac{1.273z^{6}}{z^{7}-0.8116z^{6}-1.059z^{5}+1.324z^{4}-0.2336z^{3}-0.6179z^{2}+0.5881z-0.1622} \quad (24)$$

$$G_{50}^{80}\left(z\right)=\frac{1.202z^{6}}{z^{7}-1.025z^{6}-0.9189z^{5}+1.603z^{4}-0.4712z^{3}-0.4902z^{2}+0.4281z-0.09746} \quad (25)$$

$$G_{50}^{70}\left(z\right)=\frac{-13.1z^{6}+4.059z^{5}+8.035z^{4}+3.931z^{3}-0.9597z^{2}+3.299z-2.081}{z^{7}-0.9154z^{6}} \quad (26)$$

Unit step response of the identified discrete time higher integer order transfer function models around different operating conditions are shown in Fig. 8. By the term ‘Amplitude’ in Fig. 8 here we refer to the output of the identified system i.e. power level. It is evident from Fig. 8 that the identified dc-gains vary widely with the operating point shifting due to high nonlinearity of the reactor point kinetics. Such wide variation in the local-linear dc-gains make this typical nonlinear process very difficult to control with step change in command using standard controller designing techniques. Saha *et al.* [9] and Das *et al.* [8] used fractional order controllers to ensure iso-damped controllers to ensure dead-beat power tracking for the reactor. Still it is well known that iso-damped control systems may exhibit oscillations if the dc-gain of the open loop system is increased to a large extent due to the process nonlinearity or operator’s mishandling. In this paper, a hyper-damped control system design has been attempted which will not only restrict oscillations in reactor power but also doubly ensure higher level of stability which is necessary for safer reactor operation. For model validation the AIC criteria is a standard tool which is reported in Table 1. The step response validation between estimated and experimental data is already reported in [8] for the continuous time case. Here, in (19)-(26) the corresponding discrete time transfer functions to generate synthetic frequency domain data have been shown. These are used for identifying FO models.

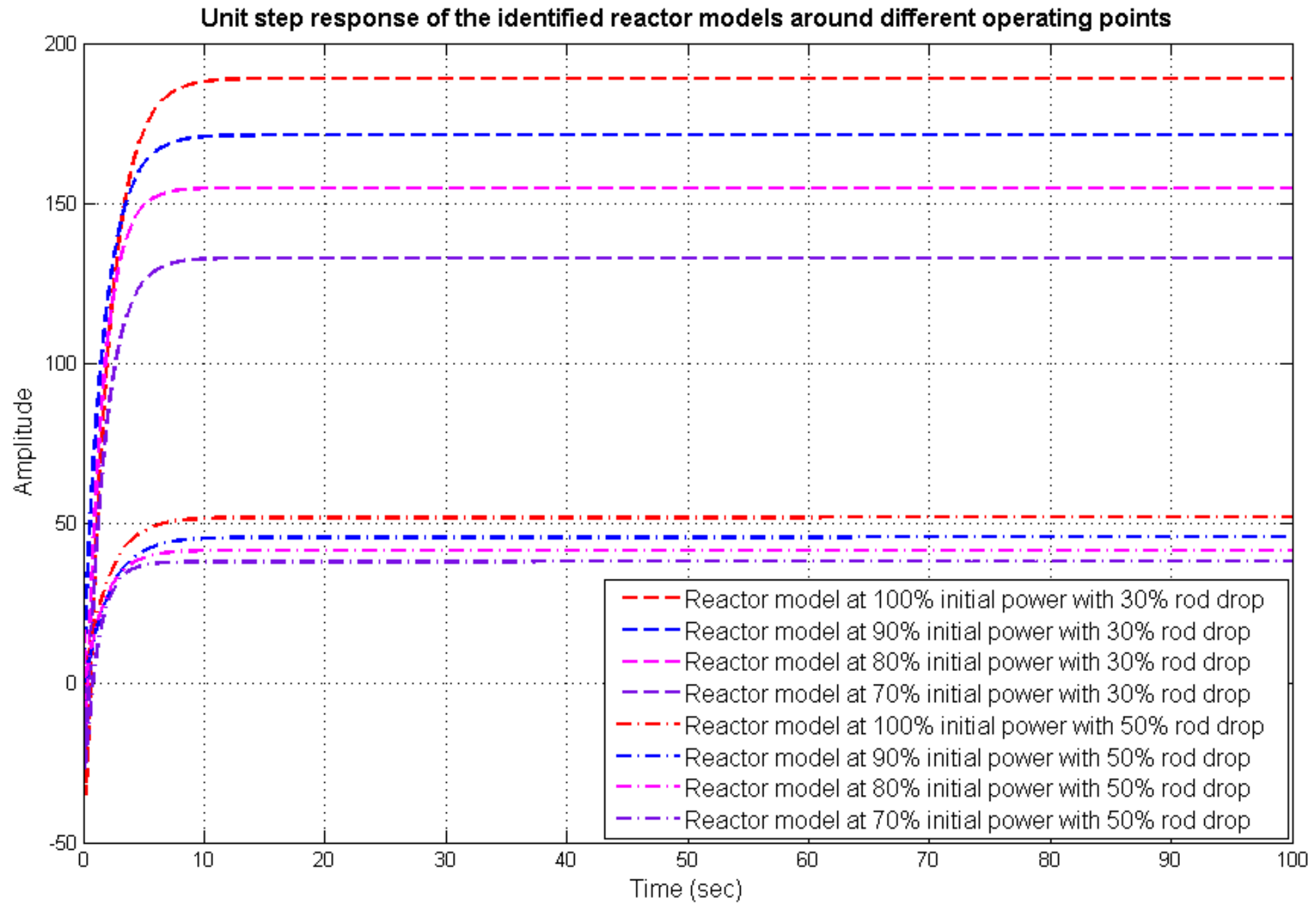

Fig. 8. Unit step response of the identified reactor models around different operating points.

## 3. Continuous order modeling of the nuclear reactor

### *3.1. Basic philosophy of continuous order system modeling*

The continuous order modeling is a new way of looking at dynamical systems which assumes that the underlying physical laws or the differential equation has a continuous distribution in its order instead of the common notion of expressing them as discrete integer order or discrete fractional order differential equations. Theoretical framework for continuous order system identification using frequency domain data has been introduced by Hartley and Lorenzo [10] to estimate all pole transfer function models and by Nazarian and Haeri [24] to estimate pole-zero models. The present approach uses the synthetic frequency domain data extracted from the identified discrete time higher integer order models (19)-(26) as suggested in [14].

There are other available techniques of fractional order system identification like FO subspace method [38], fractional orthogonal basis function [39], fractional Laguerre basis function [40], output-error technique [2], [41]-[43], frequency domain methods [44], [13]-[16] etc. Time domain system identification using linear and nonlinear estimators has been studied in [45] in the presence of fractional Gaussian noise (fGn) which shows nonlinear Hammerstein-Wiener class of estimators are well equipped in accurate modeling of linear systems over nonlinear ARX and other linear LSE variants in the presence of fGn. Multiple Riemann sheet approach of fractional order system [46] and its applications in time and frequency domain continuous order system identification have been illustrated in [47]. Basic concepts of root locus for fractional order systems [48], existence of poles in multiple Riemann sheets [49] and extensions [50]-[51] are becoming increasingly popular among the contemporary research community.

While dealing with integer order modeling for real systems, a priori knowledge of the system's highest order plays pivotal role in determining the number of parameters to be estimated. Since, the nuclear reactor whose dynamics is governed by nonlinear point kinetic equations can not be treated in a same fashion. Here, order of the system is related to the approximate linearized models around each operating point corresponding to a specific rod drop level and initial reactor power as reported in (19)-(26). However, with fractional order models, the maximum number of parameters to be estimated increases drastically depending on the sampling of order distribution i.e. the commensurate order ($q$). For continuous order models which can be visualized as a fractional order model with very small commensurate order ($q \rightarrow 0$), the parameters to be estimated can take an infinitely large number since the system's orders can be thought to have a continuous distribution between zero and the highest (fractional) order. For this reason, the sampled order distribution should be finite which can be viewed like a trade-off between low commensurate order and improvement in modeling accuracy. Also, in the pioneering work on continuous order system identification it was suggested that only few of the orders with high value of the associated coefficient should be considered as dominant orders and rest of terms from the order distribution can be ignored. This concept has been modified in another way to find out the optimum orders of the reactor under step-back in [8] using two flexible order templates.

### *3.2. Continuous order modeling with Levy's method and its variants*

It was suggested by Valerio and Sa da Costa [14] that fractional order modeling from time response data can be done in the following two steps:

a) Firstly identifying an accurate discrete time higher integer order models with available time domain system identification techniques which may be used then to generate synthetic frequency domain data.
b) This frequency domain data can now be used to develop fractional order models with the available algorithms like Levy's method and its variants like improvement with Vinagre's weight etc. [15]-[16] having minimum modeling error.

In the present work, the frequency response of identified discrete time models (19)-(26) has been used to identify FO models since the frequency domain techniques are popular for estimating FO models. The basics of fractional order system identification with Levy's method and its variants [13], [15] are briefly introduced in the next subsection.

#### *3.2.1. Levy's identification method for fractional order systems [13], [15]*

Assuming a linear system is described by a transfer function $G$ having a frequency response $G(j\omega)$, the identification consists of finding out another transfer function of the form

$$\widehat{G}(s) = \frac{b_0 + b_1 s^q + b_2 s^{2q} + \cdots + b_m s^{mq}}{a_0 + a_1 s^q + a_2 s^{2q} + \cdots + a_n s^{nq}} = \frac{\sum_{k=0}^{m} b_k s^{kq}}{\sum_{k=0}^{n} a_k s^{kq}} \tag{27}$$

where, the orders $m$ and $n$ of the numerator and denominator respectively are user specified and $q$ is the order of fractional derivative. Now, setting $a_0 = 1$ in (27) we have the corresponding frequency response of the identified model as

$$\widehat{G}(j\omega) = \frac{\sum_{k=0}^{m} b_k (j\omega)^{kq}}{1 + \sum_{k=1}^{n} a_k (j\omega)^{kq}} = \frac{N(j\omega)}{D(j\omega)} = \frac{\alpha(\omega) + j\beta(\omega)}{\sigma(\omega) + j\tau(\omega)} \tag{28}$$

where, $N$ and $D$ are complex valued and $\alpha, \beta, \sigma, \tau$ are real valued. From (28) we have

$$\begin{aligned}
\alpha(\omega) &= \sum_{k=0}^{m} b_k \operatorname{Re}\left[(j\omega)^{kq}\right] \\
\sigma(\omega) &= \sum_{k=0}^{n} a_k \operatorname{Re}\left[(j\omega)^{kq}\right] = 1 + \sum_{k=1}^{n} a_k \operatorname{Re}\left[(j\omega)^{kq}\right] \\
\beta(\omega) &= \sum_{k=0}^{m} b_k \operatorname{Im}\left[(j\omega)^{kq}\right] \\
\tau(\omega) &= \sum_{k=0}^{n} a_k \operatorname{Im}\left[(j\omega)^{kq}\right] = \sum_{k=1}^{n} a_k \operatorname{Im}\left[(j\omega)^{kq}\right]
\end{aligned} \tag{29}$$

The error between the identified model and the actual system is then given by

$$\varepsilon(j\omega) = G(j\omega) - \frac{N(j\omega)}{D(j\omega)} \tag{30}$$

Since it is difficult to choose parameters in (27) such that the error in (30) is minimized, Levy's method minimizes the square of the norm of

$$E(j\omega) := \varepsilon(j\omega) D(j\omega) = G(j\omega) D(j\omega) - N(j\omega) \tag{31}$$

which gives a set of normal equations having a simpler solution method.
Dropping the frequency argument $\omega$ to obtain a simple notation of (31) we get

$$\begin{aligned}
E &= GD - N \\
&= \left[\operatorname{Re}(G) + j \operatorname{Im}(G)\right](\sigma + j\tau) - (\alpha + j\beta) \\
&= \left[\operatorname{Re}(G)\sigma - \operatorname{Im}(G)\tau - \alpha\right] + j\left[\operatorname{Re}(G)\tau + \operatorname{Im}(G)\sigma - \beta\right]
\end{aligned} \tag{32}$$

Also, $|E|^2 = \left[\operatorname{Re}(G)\sigma - \operatorname{Im}(G)\tau - \alpha\right]^2 + \left[\operatorname{Re}(G)\tau + \operatorname{Im}(G)\sigma - \beta\right]^2$ (33)

Now, differentiating (33) with respect to one of the coefficients $b_k, k \in \{0,1,...,m\}$ or $a_k, k \in \{0,1,...,n\}$, and putting the derivative as zero, we have $\frac{\partial |E|^2}{\partial b_k} = 0$ i.e.

$$\left[\operatorname{Re}(G)\sigma - \operatorname{Im}(G)\tau - \alpha\right]\operatorname{Re}\left[(j\omega)^{kq}\right] + \left[\operatorname{Re}(G)\tau + \operatorname{Im}(G)\sigma - \beta\right]\operatorname{Im}\left[(j\omega)^{kq}\right] = 0 \tag{34}$$

Similarly, $\frac{\partial |E|^2}{\partial a_k} = 0$ yields

$$\sigma\left\{\left[\mathrm{Im}(G)\right]^2+\left[\mathrm{Re}(G)\right]^2\right\}\mathrm{Re}\left[(j\omega)^{kq}\right]+\tau\left\{\left[\mathrm{Im}(G)\right]^2+\left[\mathrm{Re}(G)\right]^2\right\}\mathrm{Im}\left[(j\omega)^{kq}\right]$$
$$+\alpha\left\{\mathrm{Im}(G)\,\mathrm{Im}\left[(j\omega)^{kq}\right]-\mathrm{Re}(G)\,\mathrm{Re}\left[(j\omega)^{kq}\right]\right\} \tag{35}$$
$$+\beta\left\{-\mathrm{Im}(G)\,\mathrm{Re}\left[(j\omega)^{kq}\right]-\mathrm{Re}(G)\,\mathrm{Im}\left[(j\omega)^{kq}\right]\right\}=0$$

The $m+1$ equations obtained from (34) and $n$ equations obtained from (35) form a linear system which can be solved to find the coefficients of (27). i.e.

$$\begin{bmatrix} A & B \\ C & D \end{bmatrix}\begin{bmatrix} b \\ a \end{bmatrix}=\begin{bmatrix} e \\ g \end{bmatrix} \tag{36}$$

where, the parameters of (36) and the corresponding expressions are detailed in [13], [15].

*3.2.2. Managing multiple frequencies*

Theoretically speaking, data from one frequency is sufficient to find a model. But in practice due to noise and other measurement inaccuracies, it is desirable to know the frequency response of the plant at more than one frequency to obtain a good identified model. There are two different approaches to deal with data from $f$ frequencies. The first approach is to sum the systems for each frequency. In this case the matrices $A, B, C, D$ and the vectors $e$ and $g$ in (36) is replaced by

$$\begin{aligned} &\widetilde{A}=\sum_{p=1}^{f}A_p\ , \quad \widetilde{B}=\sum_{p=1}^{f}B_p\ , \quad \widetilde{C}=\sum_{p=1}^{f}C_p\ , \\ &\widetilde{D}=\sum_{p=1}^{f}D_p\ , \quad \tilde{e}=\sum_{p=1}^{f}e_p\ , \quad \tilde{g}=\sum_{p=1}^{f}g_p\ . \end{aligned} \tag{37}$$

where, $A_p, B_p, C_p, D_p, e_p, g_p$ are given by (36) for a particular frequency $\omega_p$. i.e., $A_p := A(\omega_p)$ and others follow similarly. The second way is to stack several systems to obtain an over-defined system. The pseudo-inverse ($[\cdot]^+$) can be used to obtain a solution to this. Thus equation (36) becomes

$$\begin{bmatrix} A_1 & B_1 \\ C_1 & D_1 \\ A_2 & B_2 \\ C_2 & D_2 \\ \vdots & \vdots \\ A_f & B_f \\ C_f & D_f \end{bmatrix}\begin{bmatrix} b \\ a \end{bmatrix}=\begin{bmatrix} e_1 \\ g_1 \\ e_2 \\ g_2 \\ \vdots \\ e_f \\ g_f \end{bmatrix} \Rightarrow \begin{bmatrix} b \\ a \end{bmatrix}=\begin{bmatrix} A_1 & B_1 \\ C_1 & D_1 \\ A_2 & B_2 \\ C_2 & D_2 \\ \vdots & \vdots \\ A_f & B_f \\ C_f & D_f \end{bmatrix}^+\begin{bmatrix} e_1 \\ g_1 \\ e_2 \\ g_2 \\ \vdots \\ e_f \\ g_f \end{bmatrix} \tag{38}$$

*3.2.3. Adaptation of Levy's algorithm using weights*

The identification method can be enhanced using weights for each of the $f$ frequencies. Then equation (37) can be modified with weights to obtain

$$\widetilde{A}=\sum_{p=1}^{f}w_pA_p\ ,\quad \widetilde{B}=\sum_{p=1}^{f}w_pB_p\ ,\quad \widetilde{C}=\sum_{p=1}^{f}w_pC_p\ ,$$
$$\widetilde{D}=\sum_{p=1}^{f}w_pD_p\ ,\quad \tilde{e}=\sum_{p=1}^{f}w_pe_p\ ,\quad \tilde{g}=\sum_{p=1}^{f}w_pg_p\ . \tag{39}$$

*3.2.4. Vinagre's method*

Levy's method has a bias and as such often results in models which have a good fit in the high frequency data, but a poor fit in the low frequency data. Weights that decrease with frequency can be used to balance this. One reasonable value of weight is given in (40) as suggested in [13].

$$w_p=\begin{cases}\dfrac{\omega_2-\omega_1}{2\omega_1^{\,2}} & \text{if}\,p=1\\ \dfrac{\omega_{p+1}-\omega_{p-1}}{2\omega_p^{\,2}} & \text{if }\ 1<p<f\\ \dfrac{\omega_f-\omega_{f-1}}{2\omega_f^{\,2}} & \text{if}\,p=f\end{cases} \tag{40}$$

Briefly, it can be stated that Vinagre's method minimizes the norm of (41) whereas Levy's method minimizes the norm of (31).

$$E'=wG(j\omega)D(j\omega)-N(j\omega) \tag{41}$$

The accuracy of the estimated models in the original Levy's method and with Vinagre's weight can be calculated as

$$J=\frac{1}{n_\omega}\sum_{i=1}^{n_\omega}\left|G(j\omega)-\widehat{G}(j\omega)\right|^2 \tag{42}$$

It is to be noted that in this application the frequency weighted version of the algorithm is preferable. In a conventional sense, for nuclear power plant controls accurate modeling only in the low frequency regions may be sufficient. But for hyper-damped control design losing small information in the farthest parts of the negative $s$-plane may also be dangerous, since in the transformed $w$-plane the whole semi-infinite negative half $s$-plane will be squeezed within a cone and the controller heavily relies on accurate modeling of the plant. That is the reason why fitting the frequency response characteristics over a wide spectrum is important.

### *3.3. Continuous time continuous order (CTCO) modeling results for the PHWR under step-back*

The present approach uses the Box-Jenkins estimator (15) based discrete time higher integer order transfer function models (19)-(26) and their frequency domain information (variation in gain and phase with frequency) to estimate continuous time fractional order models with Levy's algorithm [15]. The commensurate order is then gradually decreased from 1.0 to 0.01 to obtain the best suited order distribution for the reactor models. The present study firstly reports the modeling accuracies for different commensurate order $q$ for the 30% rod drop models (Table 2) and 50% rod drop models (Table 3). As representative cases, the continuous order distributions i.e. variation in numerator and denominator coefficients for the identified models of the form (27) with the sampled orders $q = \{0.25, 0.1, 0.01\}$ are shown in Fig. 9-11. Also it is interesting to note that the results reported in Fig. 9-11 widely differ with the state of the art techniques in continuous order system identification [10], [24]. Hartley and Lorenzo [10] reported many closed form solutions for continuous order transfer functions for ideal order distribution curves. In [10], it is assumed that the order distribution of numerator and denominator of any FO transfer function representing a stable physical system is always positive so that a smooth curve can be fitted through those discrete sampled order distributions. The pioneering work [10] has mostly given closed form transfer function representation in terms of *r*-Laplace transform (logarithmic Laplace transform) for ideal order distributions like Gaussian/triangular around a dominant order, uniform, saw-tooth or spiky etc. The present study investigates the rationale behind assuming such structured order distribution curves for a practical system i.e. a nuclear reactor under step-back.

It is observed from Fig. 9-11 that for varying level of sampling in the order distribution i.e. commensurate order $q$, the distribution of the coefficients associated with the numerator/denominator of continuous order model is not at all smooth so that their variation can be approximated with available curve fitting techniques. In fact, the coefficients $K(q)$ widely varies in both positive and negative direction still representing a stable transfer function model. This is justified due to the fact that the stable poles in $s$-plane gets mapped onto the $w$-plane, associated with the corresponding commensurate order using the relation $w = s^q$ [27], [1]. Therefore for fractional order transfer function models negative terms in the denominator do not always represent unstable dynamics.

For instance a characteristic polynomial of type say $\left(s + a\sqrt{s} + 1\right)$ i.e. a first order system with half order element, is stable for $a \in (0, -\sqrt{2})$, that is having negative constants in indicial polynomial. For these values of $\left(-\sqrt{2} < a < 0\right)$, the system is stable ultra-damped. While $a \in (0, 2)$ this fractional order system is stable with roots in secondary Riemann sheet as a hyper-damped system. While still the system is stable with roots in secondary Riemann sheet, as ultra-damped system when $a > 2$ . This is one example contrary to integer order systems, when the constants need to be always positive for stability. Thus for fractional order systems the constants if indicial polynomial can be negative still giving stability. In this example, while $a < -\sqrt{2}$ , the system is unstable.

Table 2
Frequency domain continuous order modeling results for 30% rod drop models

| Model | Commensurate order (q) | Accuracy of Identification algorithms (*J*) | |
|---|---|---|---|
| | | Levy | Levy with Vinagre's weight |
| $G_{30}^{100}$ | 1.0 | $2.4777\times10^{8}$ | $2.0043\times10^{10}$ |
| | 0.5 | 1.0848 | 0.3349 |
| | 0.25 | 0.1164 | $3.6721\times10^{-6}$ |
| | 0.1 | $2.1373\times10^{-5}$ | $3.0706\times10^{-6}$ |
| | 0.05 | $2.1571\times10^{-6}$ | $2.3511\times10^{-6}$ |
| | 0.02 | $1.1557\times10^{-5}$ | $7.7031\times10^{-6}$ |
| | 0.01 | $1.4561\times10^{-5}$ | $4.4239\times10^{-6}$ |
| $G_{30}^{90}$ | 1.0 | $5.5254\times10^{7}$ | $4.2096\times10^{8}$ |
| | 0.5 | 13.016 | 9.0693 |
| | 0.25 | 0.1577 | 0.013 |
| | 0.1 | $6.8638\times10^{-4}$ | 0.0024 |
| | 0.05 | 0.0018 | 0.0044 |
| | 0.02 | 0.0247 | 0.0177 |
| | 0.01 | 0.002 | 0.0027 |
| $G_{30}^{80}$ | 1.0 | $2.1284\times10^{7}$ | $2.9793\times10^{9}$ |
| | 0.5 | 49.6226 | 25.3173 |
| | 0.25 | $5.2691\times10^{-4}$ | $8.0972\times10^{-4}$ |
| | 0.1 | 0.0014 | $9.3118\times10^{-4}$ |
| | 0.05 | $5.8145\times10^{-4}$ | $1.4332\times10^{-4}$ |
| | 0.02 | 0.0043 | $3.846\times10^{-6}$ |
| | 0.01 | $4.7787\times10^{-4}$ | $3.1073\times10^{-4}$ |
| $G_{30}^{70}$ | 1.0 | 2.1123×107 | $4.6137\times10^{9}$ |
| | 0.5 | 28.1039 | 1.2883 |
| | 0.25 | 0.0011 | 0.0065 |
| | 0.1 | $2.4398\times10^{-4}$ | 0.0014 |
| | 0.05 | $3.236\times10^{-4}$ | $7.1466\times10^{-4}$ |
| | 0.02 | $2.134\times10^{-4}$ | $5.9372\times10^{-4}$ |
| | 0.01 | $2.1278\times10^{-4}$ | 0.0044 |

Also, the concept of dominant order, introduced in [10] i.e. retention of high values in the order distribution while neglecting small coefficients might not always lead to stability or preservation of the original dynamics. Therefore, we find that the notion of giving more importance to numerically large values in the order distribution and neglecting small coefficients [10], [52] is not always correct. Instead, an equivalent compressed model should be searched for using an optimization some technique that optimally represent all the dynamics associated with individual sampled orders and their weights (coefficients) into a compact template while keeping the order of the derivatives flexible [12], [8].

Table 3
Frequency domain continuous order modeling results for 50% rod drop models

| **Model** | **Commensurate order (q)** | **Accuracy of Identification algorithms (*J*)** | |
|---|---|---|---|
| | | **Levy** | **Levy with Vinagre's weight** |
| $G_{50}^{100}$ | 1.0 | $5.1338\times10^{6}$ | $1.2883\times10^{7}$ |
| | 0.5 | 0.8169 | 0.0818 |
| | 0.25 | $7.5497\times10^{-4}$ | 0.0015 |
| | 0.1 | $1.5193\times10^{-5}$ | $2.1676\times10^{-4}$ |
| | 0.05 | $4.5358\times10^{-5}$ | $7.5941\times10^{-4}$ |
| | 0.02 | $8.9297\times10^{-6}$ | $9.8821\times10^{-4}$ |
| | 0.01 | $6.7569\times10^{-6}$ | $3.6454\times10^{-5}$ |
| $G_{50}^{90}$ | 1.0 | $1.7614\times10^{5}$ | $2.2024\times10^{7}$ |
| | 0.5 | 0.2629 | 1.3884 |
| | 0.25 | 0.0055 | 0.0042 |
| | 0.1 | 0.0024 | 0.0046 |
| | 0.05 | 0.0016 | 0.0041 |
| | 0.02 | $6.2712\times10^{-5}$ | 0.0011 |
| | 0.01 | 0.0015 | 0.0044 |
| $G_{50}^{80}$ | 1.0 | $3.6410\times10^{5}$ | $2.7793\times10^{7}$ |
| | 0.5 | 0.4897 | 0.6018 |
| | 0.25 | 8.7292 | 0.2087 |
| | 0.1 | $9.9278\times10^{-4}$ | 0.0043 |
| | 0.05 | 0.0018 | 0.003 |
| | 0.02 | 0.0056 | 0.0086 |
| | 0.01 | $7.304\times10^{-4}$ | 0.0059 |
| $G_{50}^{70}$ | 1.0 | $3.3248\times10^{6}$ | $4.5481\times10^{8}$ |
| | 0.5 | 19.5121 | 18.0993 |
| | 0.25 | 0.001 | 0.0036 |
| | 0.1 | $3.329\times10^{-4}$ | $5.8482\times10^{-4}$ |
| | 0.05 | $1.3247\times10^{-4}$ | 0.0196 |
| | 0.02 | $6.1855\times10^{-4}$ | 0.0017 |
| | 0.01 | 0.002 | $3.2374\times10^{-4}$ |

It is also found that for very low value of the commensurate orders the system matrices which needs to be inverted within the algorithm, become close to singular due to their drastic increase in size. As a result estimation problem becomes more and more inconsistent with increase in computational complexity. As a trade-off between better accuracy and low complexity of the model we have restricted the commensurate order as $q = 0.25$ and the corresponding FO reactor models are reported in (43)-(50). From Table 2 and 3, it is also evident that the frequency domain identification accuracy of the continuous order models increases if the commensurate order ($q$) is decreased gradually, so that the whole sampled order distribution can be seen in a finer resolution. But for $q < 0.1$ the argument of fall in accuracy with finer resolution becomes inconsistent due to the fact that the system matrices become larger and this also increases the parametric variance of the estimated model coefficients. Fig. 12 shows that the frequency domain validation of the identified discrete time higher integer order reactor models (19)-(26) with the continuous time continuous order reactor models (43)-(50) considering a commensurate order of $q = 0.25$, as discussed earlier. The Bode diagram in Fig. 12 shows that the CTCO models have efficiently described the frequency domain information of the discrete time integer order models up to the corresponding Nyquist frequency. It is also interesting to note from the continuous order reactor models in (43)-(50) that in the presence of other fractional order elements a stable system can have the highest fractional order more than two as reported in Das *et al.* [12], [8] and has been assumed here as $mq = nq = 2.5$ in equation (27) to estimate continuous order models of the reactor under step-back.

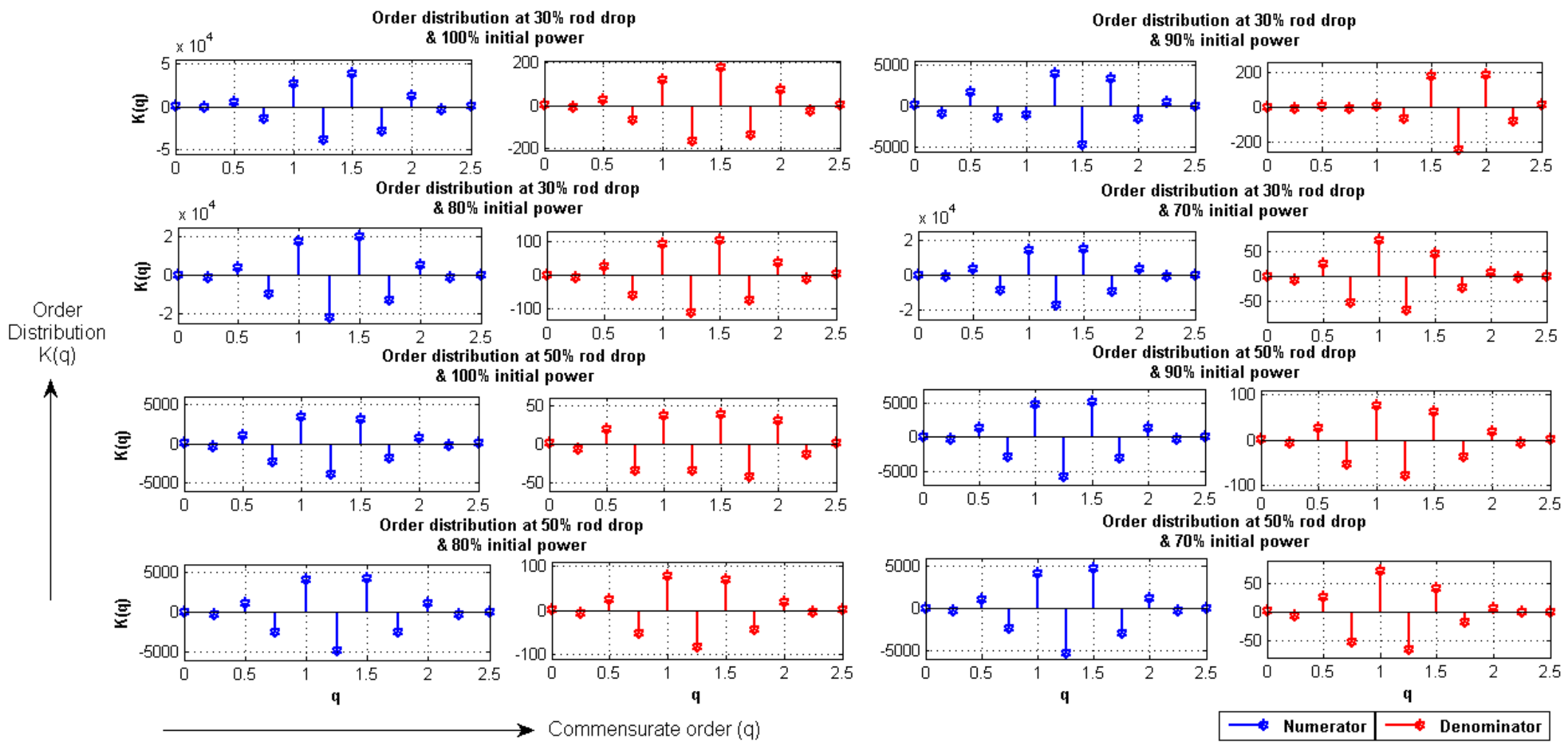

Fig. 9. Order distribution of the identified models having commensurate order q=0.25.

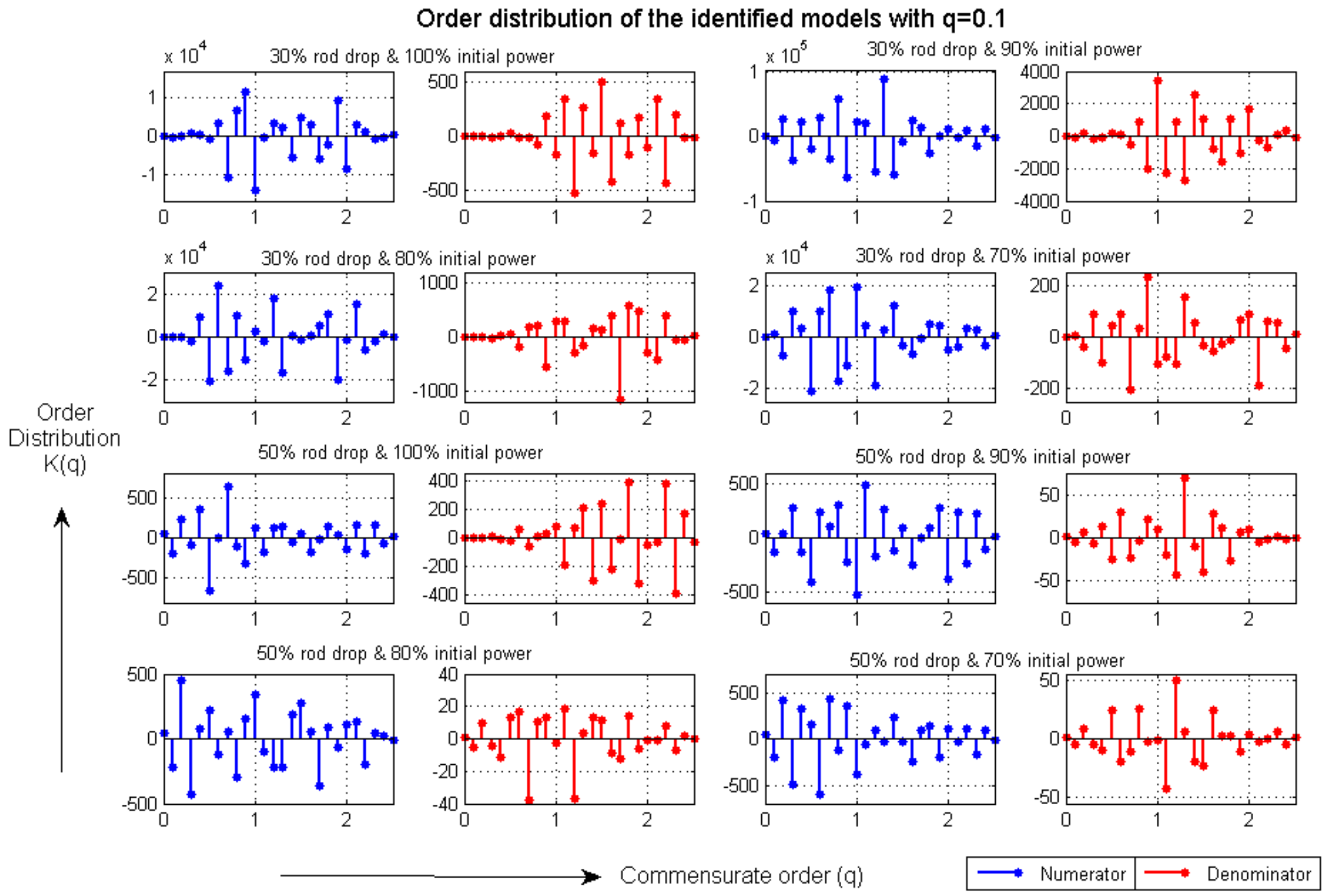

Fig. 10. Order distribution of the identified models having commensurate order q=0.1.

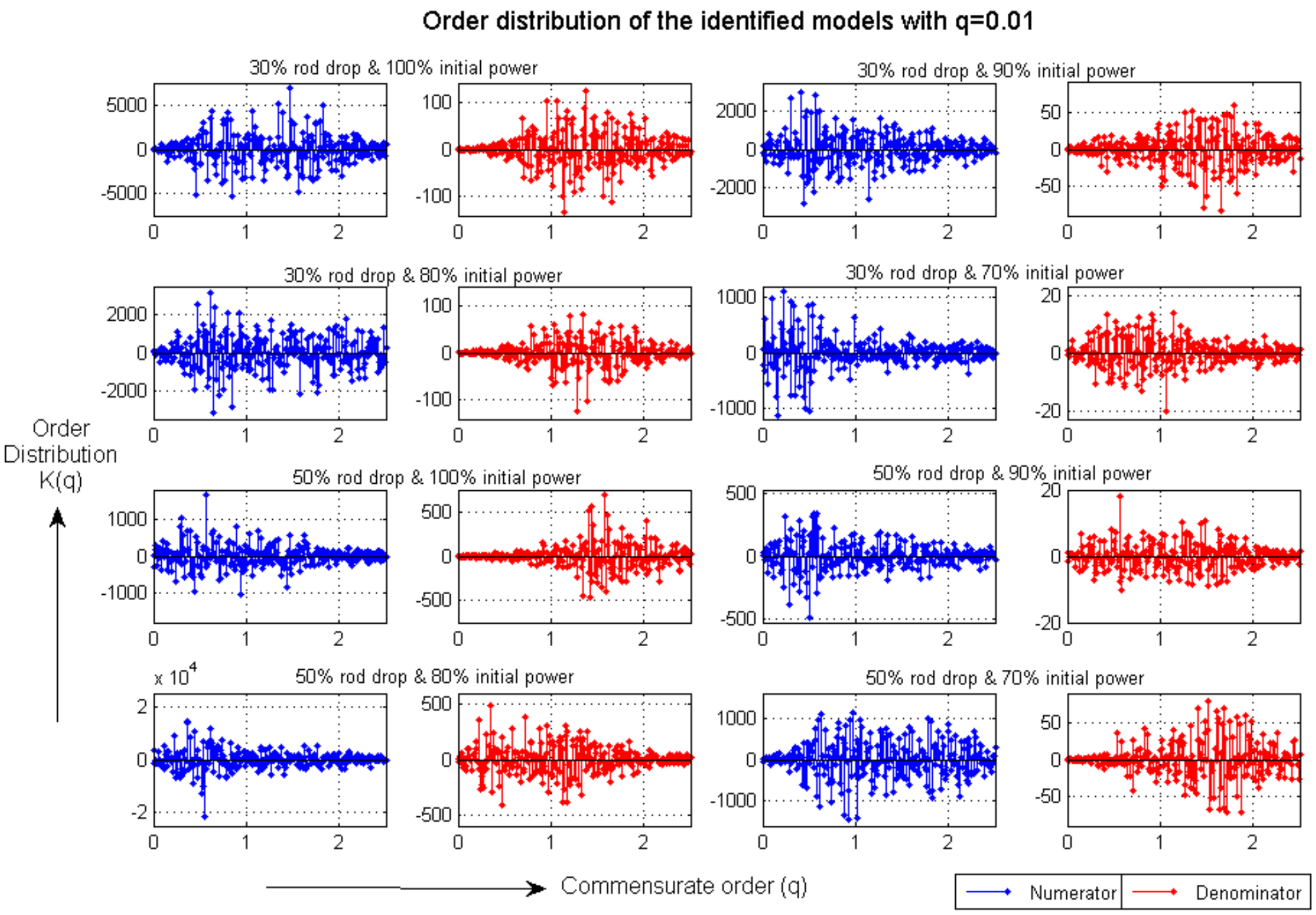

Fig. 11. Order distribution of the identified models having commensurate order q=0.01.

$$\tilde{G}_{30}^{100}(s)=\frac{\begin{array}{c}442.8093s^{2.5}-3584.6003s^{2.25}+12929.3346s^{2}-27507.7124s^{1.75}+38563.6082s^{1.5}\\-37756.4006s^{1.25}+26716.4613s-13862.9884s^{0.75}+5205.5968s^{0.5}-1343.1001s^{0.25}+189.2362\end{array}}{\begin{array}{c}4.0473s^{2.5}-25.1112s^{2.25}+74.4725s^{2}-136.868s^{1.75}+175.4011s^{1.5}-166.332s^{1.25}\\+120.5556s-66.2508s^{0.75}+26.7975s^{0.5}-7.0595s^{0.25}+1\end{array}} \quad (43)$$

$$\tilde{G}_{30}^{90}(s)=\frac{\begin{array}{c}-51.2735s^{2.5}+412.3702s^{2.25}-1527.1234s^{2}+3359.9726s^{1.75}-4668.634s^{1.5}+3848.8028s^{1.25}\\-1100.2029s-1298.9864s^{0.75}+1707.7228s^{0.5}-852.8921s^{0.25}+171.3044\end{array}}{\begin{array}{c}14.2649s^{2.5}-78.4487s^{2.25}+186.3603s^{2}-241.2807s^{1.75}+175.715s^{1.5}-64.2108s^{1.25}+8.7947s\\-6.8852s^{0.75}+9.8484s^{0.5}-4.9694s^{0.25}+1\end{array}} \quad (44)$$

$$\tilde{G}_{30}^{80}(s)=\frac{\begin{array}{c}149.8262s^{2.5}-1337.0308s^{2.25}+5358.5713s^{2}-12740.6685s^{1.75}+20018.8426s^{1.5}-21943.9262s^{1.25}\\+17283.403s-9902.8347s^{0.75}+4073.4683s^{0.5}-1115.0882s^{0.25}+154.9787\end{array}}{\begin{array}{c}1.8404s^{2.5}-12.1146s^{2.25}+37.607s^{2}-73.2914s^{1.75}+102.1556s^{1.5}-109.157s^{1.25}\\+91.0873s-57.4605s^{0.75}+25.5829s^{0.5}-7.1556s^{0.25}+1\end{array}} \quad (45)$$

$$\widetilde{G}_{30}^{70}(s)=\frac{\begin{array}{c}89.9109s^{2.5}-846.0716s^{2.25}+3584.1383s^{2}-9003.7538s^{1.75}+14897.4822s^{1.5}-17095.4976s^{1.25}\\+13995.0615s-8277.9642s^{0.75}+3492.1868s^{0.5}-968.8845s^{0.25}+133.1494\end{array}}{\begin{array}{c}0.16026s^{2.5}-1.4727s^{2.25}+6.8837s^{2}-20.9657s^{1.75}+44.6741s^{1.5}-67.554s^{1.25}\\+71.8732s-52.4714s^{0.75}+25.1083s^{0.5}-7.2106s^{0.25}+1\end{array}} \tag{46}$$

$$\widetilde{G}_{50}^{100}(s)=\frac{\begin{array}{c}18.416s^{2.5}-171.2393s^{2.25}+724.5365s^{2}-1843.103s^{1.75}+3145.2927s^{1.5}-3813.3739s^{1.25}\\+3393.7524s-2241.2139s^{0.75}+1070.7003s^{0.5}-335.5444s^{0.25}+51.8497\end{array}}{\begin{array}{c}2.2383s^{2.5}-12.562s^{2.25}+30.4049s^{2}-41.347s^{1.75}+38.0912s^{1.5}-34.1689s^{1.25}\\+36.7103s-33.2323s^{0.75}+19.3428s^{0.5}-6.3842s^{0.25}+1\end{array}} \tag{47}$$

$$\widetilde{G}_{50}^{90}(s)=\frac{\begin{array}{c}35.2472s^{2.5}-315.5662s^{2.25}+1284.5154s^{2}-3133.4332s^{1.75}+5090.714s^{1.5}-5798.5053s^{1.25}\\+4750.6803s-2818.5909s^{0.75}+1186.9863s^{0.5}-327.4343s^{0.25}+45.4763\end{array}}{\begin{array}{c}0.99301s^{2.5}-5.9022s^{2.25}+17.8417s^{2}-37.2445s^{1.75}+60.5457s^{1.5}-77.8108s^{1.25}\\+76.1043s-53.4213s^{0.75}+25.1562s^{0.5}-7.1464s^{0.25}+1\end{array}} \tag{48}$$

$$\widetilde{G}_{50}^{80}(s)=\frac{\begin{array}{c}26.6578s^{2.5}-244.6936s^{2.25}+10204001s^{2}-25450691s^{1.75}+4215.341s^{1.5}-4878.1213s^{1.25}\\+404884s-24322661s^{0.75}+1040.2897s^{0.5}-292.7986s^{0.25}+41.4599\end{array}}{\begin{array}{c}0.65703s^{2.5}-5.0526s^{2.25}+185352s^{2}-425015s^{1.75}+68.4999s^{1.5}-82.6806s^{1.25}\\+761678s-519011s^{0.75}+24.3433s^{0.5}-7.0195s^{0.25}+1\end{array}} \tag{49}$$

$$\widetilde{G}_{50}^{70}(s)=\frac{\begin{array}{c}31.2553s^{2.5}-290.6344s^{2.25}+1210.6045s^{2}-2973.7433s^{1.75}+4783.2853s^{1.5}-5309.0061s^{1.25}\\+4195.4509s-2407.8342s^{0.75}+999.2335s^{0.5}-276.5767s^{0.25}+37.8298\end{array}}{\begin{array}{c}0.14397s^{2.5}-1.1345s^{2.25}+5.4393s^{2}-18.1782s^{1.75}+41.8715s^{1.5}-66.0764s^{1.25}\\+71.4213s-52.2851s^{0.75}+25.0496s^{0.5}-7.2277s^{0.25}+1\end{array}} \tag{50}$$

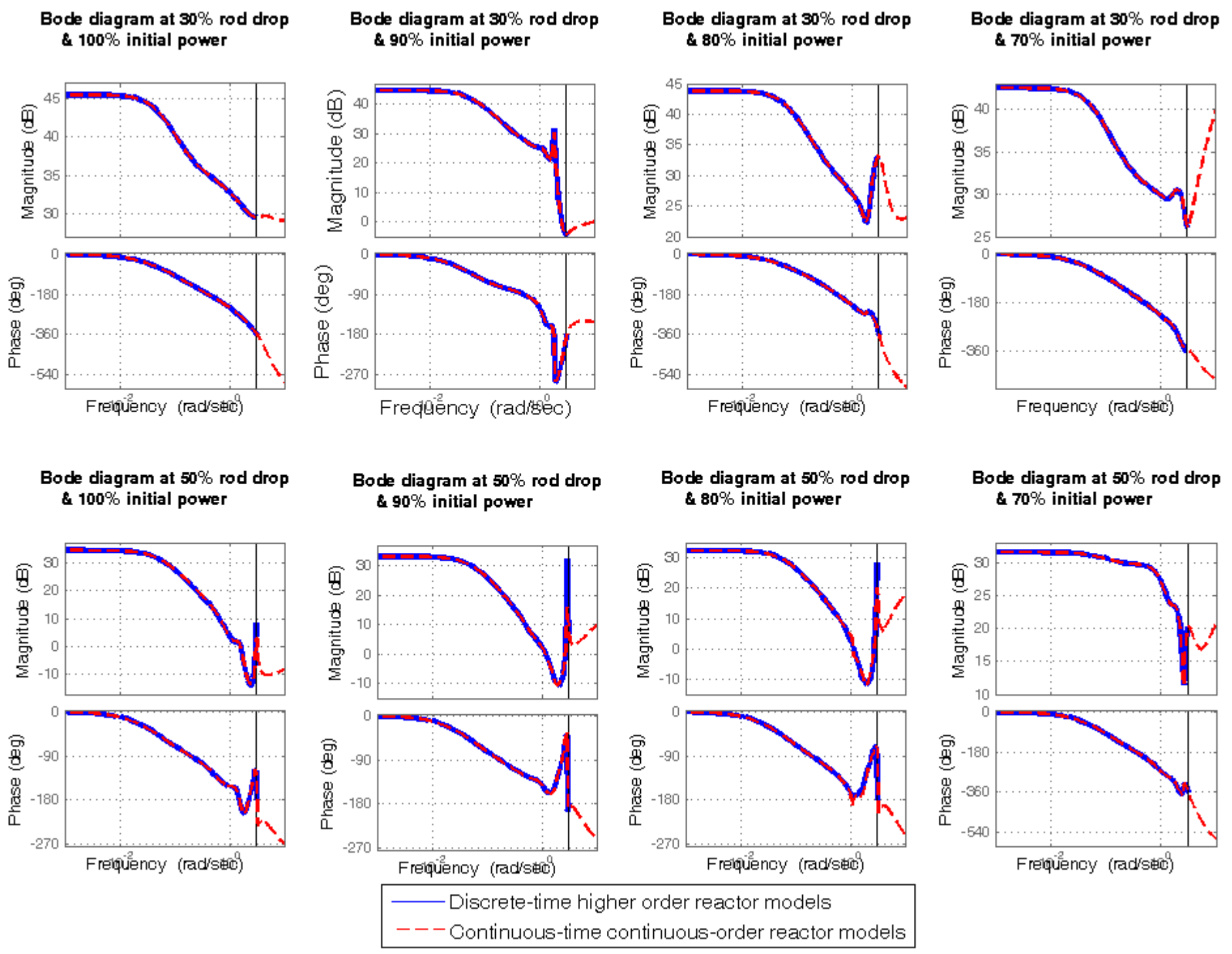

Fig. 12. Frequency domain validation for continuous order modeling.

### *3.4. Analysis of identified continuous order reactor models*

A closer look at the identified continuous order models (43)-(50) with commensurate order $q = 0.25$ is now needed, in the line of stability of those models. In order to do so, the basic concept of stability and dynamics in complex $w$-plane has been discussed first [1], [27]. Let us assume that a fractional order transfer function takes the form (27) with commensurate order $q$. If $\lambda_i$ be the poles of the FO model then the system is stable for the condition, that is $\left|\arg(\lambda_i)\right| > \pi q/2$. It has been illustrated in Fig. 13 that $\left|\arg(\lambda_i)\right| < \pi q/2$ yields unstable dynamics. With the concept of fractional order systems the higher Riemann sheets come into play i.e. poles lying in the region $\left|\arg(\lambda_i)\right| > \pi q$. These concepts can not be visualized using conventional integer order concepts of poles, zeros or root locus and therefore the corresponding fractional order dynamics and stability versions should be used [49]. FO systems with $\pi q < \left|\arg(\lambda_i)\right| < \pi$ are known as hyper-damped whereas with $\left|\arg(\lambda_i)\right| = \pi$ it will be termed as ultra-damped system.

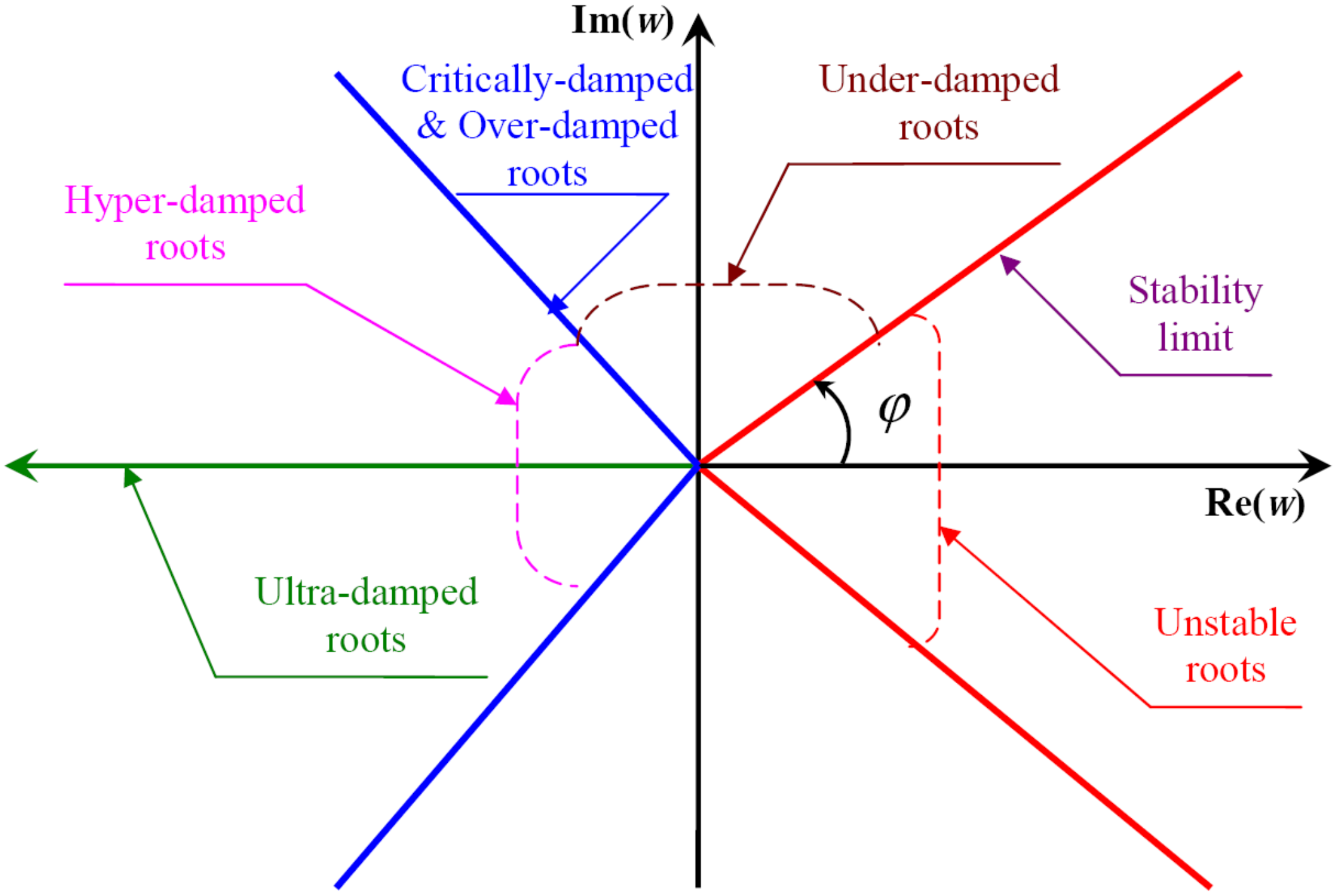

Fig. 13. Existence of hyper-damped and ultra-damped poles in higher Riemann sheets.

The concept has been visualized in Fig. 13 in a self-explanatory manner. It can be observed that the whole negative half of the $s$-plane gets compressed within the region $\pi q/2 < |\arg(\lambda_i)| < \pi q$ in complex $w$-plane. Also, some additional higher Riemann sheets have appeared with the possibility of existing poles in these zones. Since, over-damped feedback control system can still go to oscillations if the gain of the open loop system is increased in a significant manner as all the branches remains in the primary Riemann sheet for integer order dynamical systems. But for the case of a fractional order system, if it is enforced in the FO controller design stage so that all the poles lie in the higher Riemann sheets then a dead-beat response can be doubly ensured as the root locus branches lie in the higher Riemann sheets and can never go to oscillation or instability even for a large variation in loop gain.

The identified CTCO models with sampled order of $q = 0.25$, stability region becomes $|\arg(\lambda_i)| < 22.5^\circ$. For the identified open loop systems (43)-(50), the argument of the poles are reported in Table 4 and the corresponding pole-zero maps are shown in Fig. 14. Table 4 shows that all of the open loop poles lie above the stability region which is also justified in Fig. 14. In fact, few of the poles lie in the under-damped region also which may lead to poor performance at high gain. All the pole angles appeared in pairs of positive and negative sign but same absolute value, since they represent complex conjugates in the $w$-plane. Few of the data in Table 4 are closer to 22.5° implying closer to marginal stability operation with the power regulator only, but none of the poles has argument less than the stability limit of 22.5°. In the next section we have tried to design a single continuous order controller which will enforce dead-beat tracking and also not let the system to go to oscillations while handling the eight set of models (43)-(50) representing the reactor at different operating condition.

Table 4:
Argument of the poles for continuous order models at different operating point

| **Level of rod drop** | **30%** | | | | **50%** | | | |
|---|---|---|---|---|---|---|---|---|
| **Initial power** | **100%** | **90%** | **80%** | **70%** | **100%** | **90%** | **80%** | **70%** |
| argument of the poles $\arg(\lambda_i)$ in degrees | 30.7877 | 22.8461 | 25.8722 | 26.8784 | 22.6624 | 22.5402 | 22.5958 | 25.079 |
| | -30.7877 | -22.8461 | -25.8722 | -26.8784 | -22.6624 | -22.5402 | -22.5958 | -25.079 |
| | 34.0734 | 26.2987 | 27.4984 | 27.2783 | 26.0995 | 32.4109 | 23.3214 | 26.864 |
| | -34.0734 | -26.2987 | -27.4984 | -27.2783 | -26.0995 | -32.4109 | -23.3214 | -26.864 |
| | 45.0014 | 30.4573 | 37.0359 | 27.5585 | 33.2098 | 44.1681 | 38.4981 | 33.6058 |
| | -45.0014 | -30.4573 | -37.0359 | -27.5585 | -33.2098 | -44.1681 | -38.4981 | -33.6058 |
| | 53.9669 | 44.9721 | 45.1178 | 45.0566 | 44.9379 | 45.0754 | 45.1169 | 45.022 |
| | -53.9669 | -44.9721 | -45.1178 | -45.0566 | -44.9379 | -45.0754 | -45.1169 | -45.022 |
| | 87.4224 | 140.7488 | 94.2025 | 71.4097 | 127.6716 | 98.8554 | 89.3358 | 80.7856 |
| | -87.4224 | -140.7488 | -94.2025 | -71.4097 | -127.6716 | -98.8554 | -89.3358 | -80.7856 |

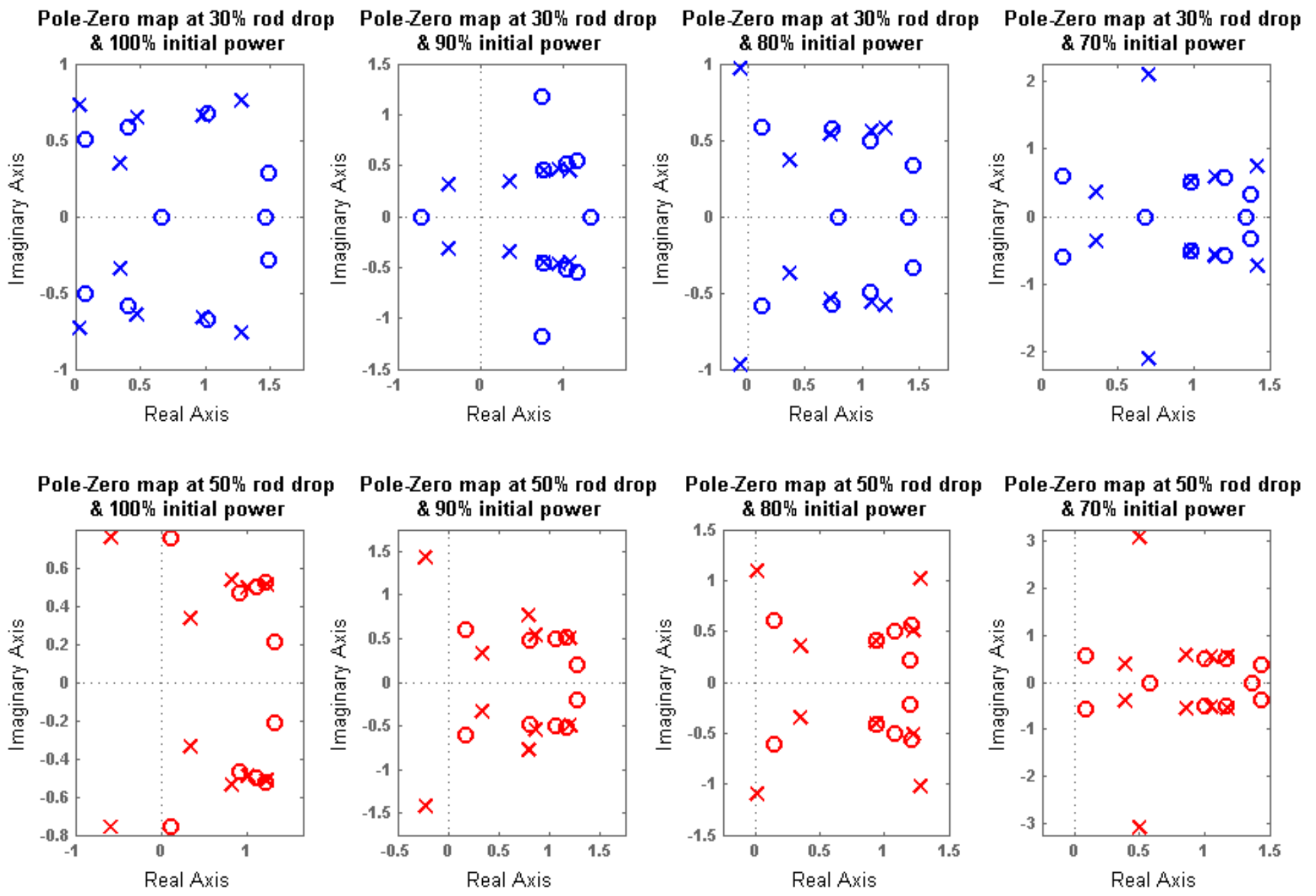

Fig. 14. Location of the identified poles and zeros of the continuous order reactor models.

## 4. Continuous order controller design for active step-back

### *4.1. Design philosophy for continuous order PID like controller*

It is already discussed that the continuous order models of the form (27) can be efficiently controlled by compensators of the structure (2) or (3) using similar controller design tasks in $w$-plane. With the structure (3) the number of controller parameters to be determined increases and set-point tracking is not guaranteed like FO lead-lag compensators [29]. Also, due to the presence of large number of FO elements in the controller structure (3) the cost of hardware realization and the complexity of its realized version will increase.

The main goal behind the design of controller in the present problem is to stabilize the dynamics of identified continuous order models in such a way so that it tracks the reference input. The tracking of the reference trajectory can be obtained by the well established methodologies that minimize the time domain performance index to find out the controller gains. Thus the optimized gains of the controller will ensure optimum time domain performance over the operating condition for which the controller is tuned. But the obtained gain will not ensure good performance or stability over the other operating points as the process gain changes with shift in operating point due to process nonlinearity. Therefore, time domain performance index optimization based FO controller design methods [12] have not been applied in the present case. Also, designing eight different controllers using the linear models at eight different test conditions and their switching is also not a feasible option as far as stability in the intermediate operating conditions are concerned. Therefore it is desirable to design a single controller which will ensure dead-beat power level tracking at all of the eight step-back conditions.

### *4.2. Continuous order controller design in an optimization frame*

In the present problem, a continuous order PID like controller of the structure (2) needs to be designed in such a way so that the poles of the closed loop systems lie outside the unstable region shown in Fig. 13. More precisely, all the closed loop poles (even at different operating point) can be pushed to the higher Riemann sheets while searching for the controller coefficients within an optimization framework. This ensures a safer reactor operation since hyper-damped poles can not exhibit oscillations even at very high gain due to nonlinearity, failure or mishandling of operator. But this extra safety feature comes at the cost of sluggishness during normal operation of the reactor as the hyper-damped poles introduce slow dynamic response.

Considering the controller structure as (2) with $q = 0.25$, an optimization based framework has been developed to search for the controller zeros while minimizing the objective function (51). The objective function (51) ensures that all the closed loop poles lie with an angle, slightly higher than $180^\circ \times q = 45^\circ$, so that increased stability due to hyper-damped poles and moderately fast time response both can be enjoyed within the same design. This makes the closed loop design faster than that with ultra-damped and hyper-damped closed loop poles which are far away from the junction between primary and secondary Riemann sheet, thus leading to very slow dynamic response. The optimization searches for controller gains (coefficients) of structure (2) until all the closed loop poles are not pushed in secondary Riemann sheet and further away i.e. in the hyper-damped zone. In (51), $\overline{\lambda}_i$ represents the $i^{th}$ closed loop pole ($i \in [1,10]$) for the eight different rod drop models (43)-(50) and norm $\|\cdot\|$ denotes the Euclidian distance.

$$\overline{J} = \left\| \left| \arg\left(\overline{\lambda}_i\right) \right| - 45^\circ \right\| \qquad (51)$$

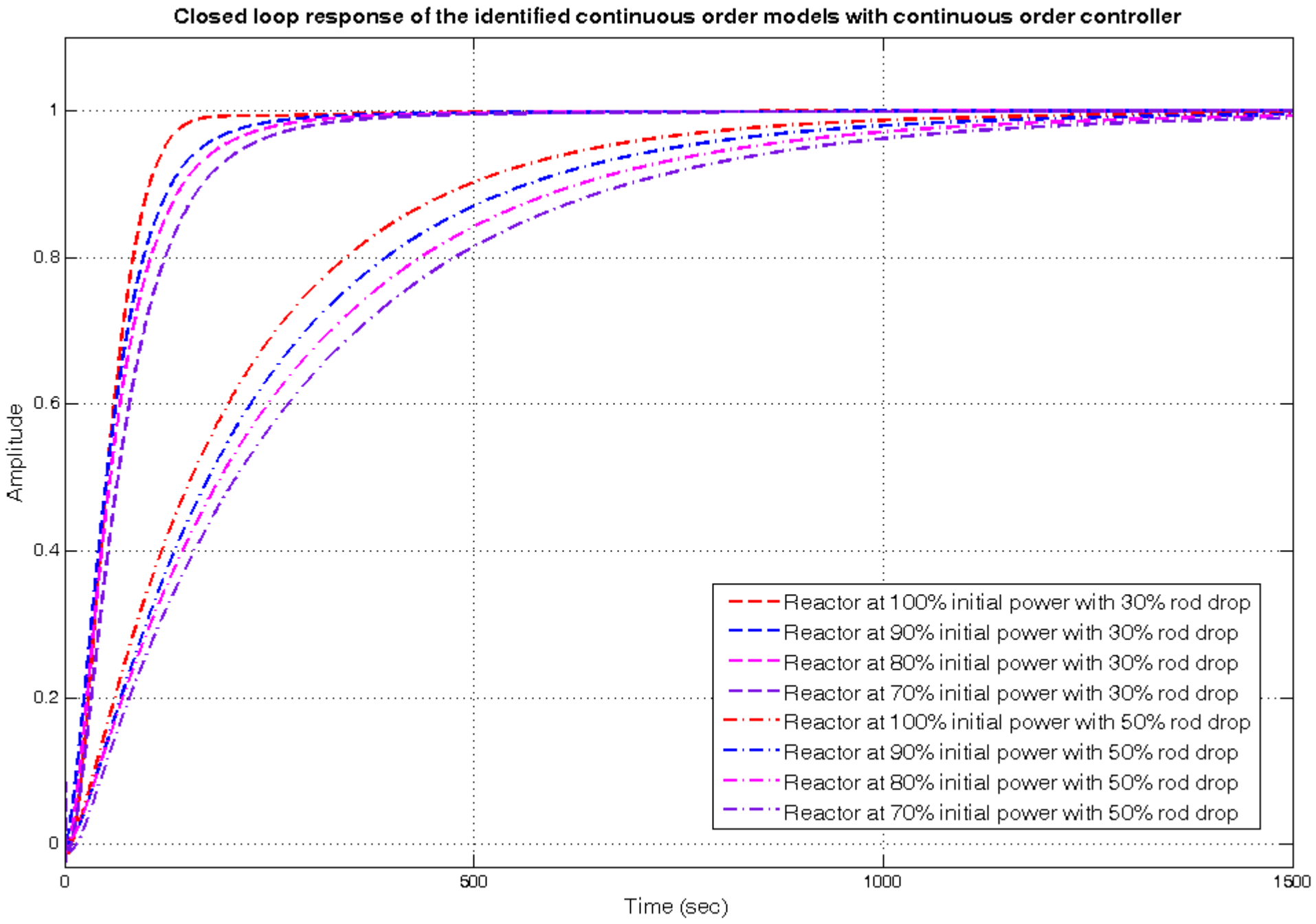

Fig. 15. Closed loop response of the identified continuous order reactor models with continuous order PID like controller.

The objective function (51) is minimized with an unconstrained Nelder-Mead simplex algorithm implemented in MATLAB's optimization toolbox [53] function *fminsearch()* with perturbed initial guess and the resulting continuous order PID like controller is reported in (52) that produces hyper-damped poles for all the eight continuous order models (43)-(50). Within the controller structure in (2) the integrator is not replaced by a fractional order one since this will lead to additional sluggishness in the system which is not desired.

$$C^{cont}(s)=\frac{\begin{bmatrix}0.5298s^{2.5}+0.2105s^{2.25}+0.9427s^{2.0}+0.6789s^{1.75}+0.4455s^{1.5}+0.0012s^{1.25}\\+0.1828s^{1.0}+0.6630s^{0.75}+0.0303s^{0.5}+0.2878s^{0.25}+0.8228\end{bmatrix}\times10^{-4}}{s} \tag{52}$$

The closed loop responses for the identified CTCO models (43)-(50) with the CTCO PID like controller (52) have been shown in Fig. 15 for unit step reference input. It is observed that the controller (52) is capable of producing dead-beat power tracking response at all operating condition though a bit sluggish time response is obtained especially at 50% rod drop conditions. Therefore, the continuous order controller (52) can be efficiently employed for the active step-back for reactor global power level control like that in Das *et al.* [8] with a $PI^{\lambda}D^{\mu}$ controller, over the present day's passive step-back mechanism. The power level tracking performances at real scale has been shown in Fig. 16 around various step-back levels and initial reactor power. It is clear that 30% drop of power can be possible within 400 seconds and also 50% drop is possible within 1600 seconds with additional safety features incorporated in the control scheme as hyper-damped closed loop poles.

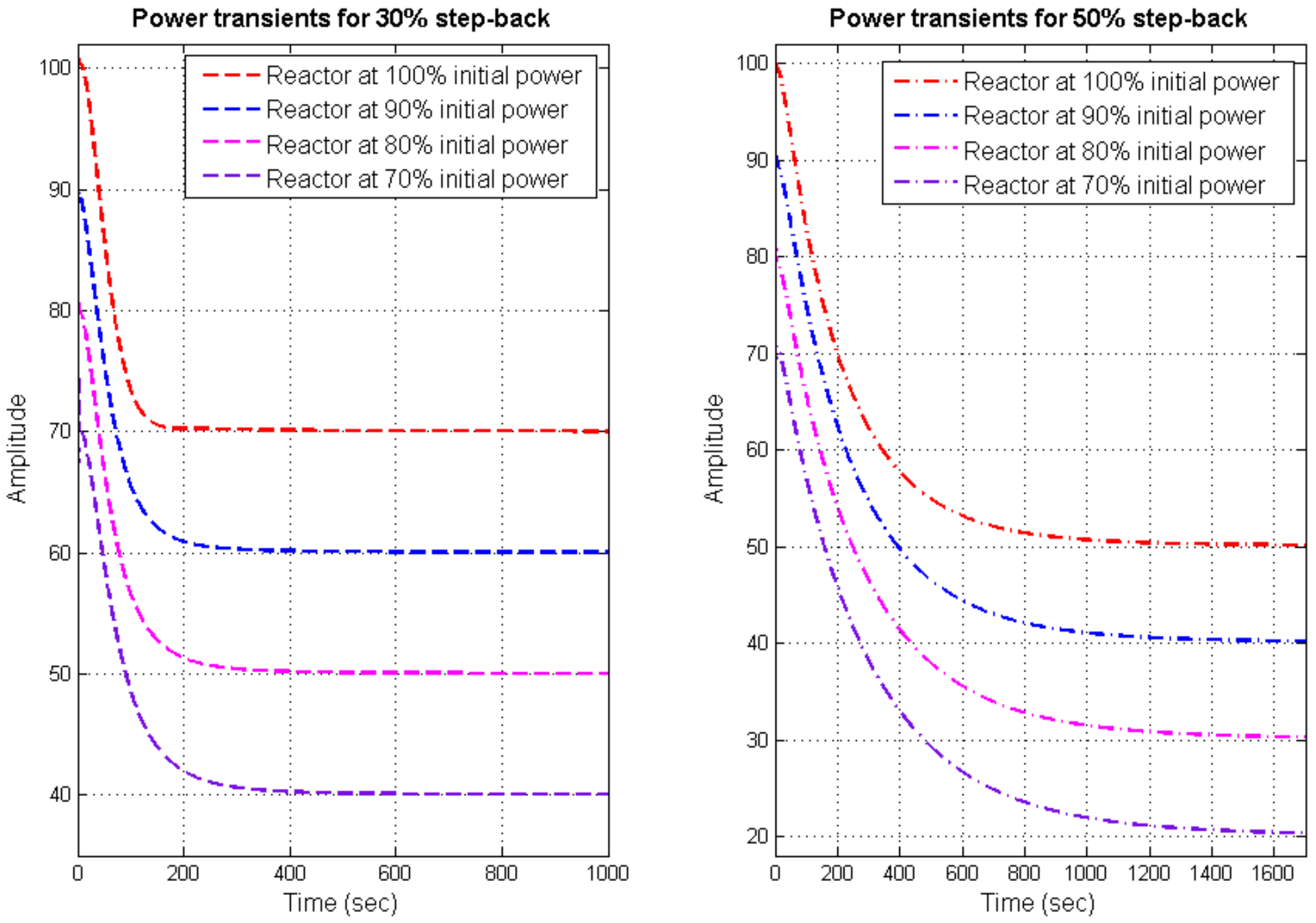

Fig. 16. 30% and 50% step-back responses of the reactor with the continuous order PID like controller.

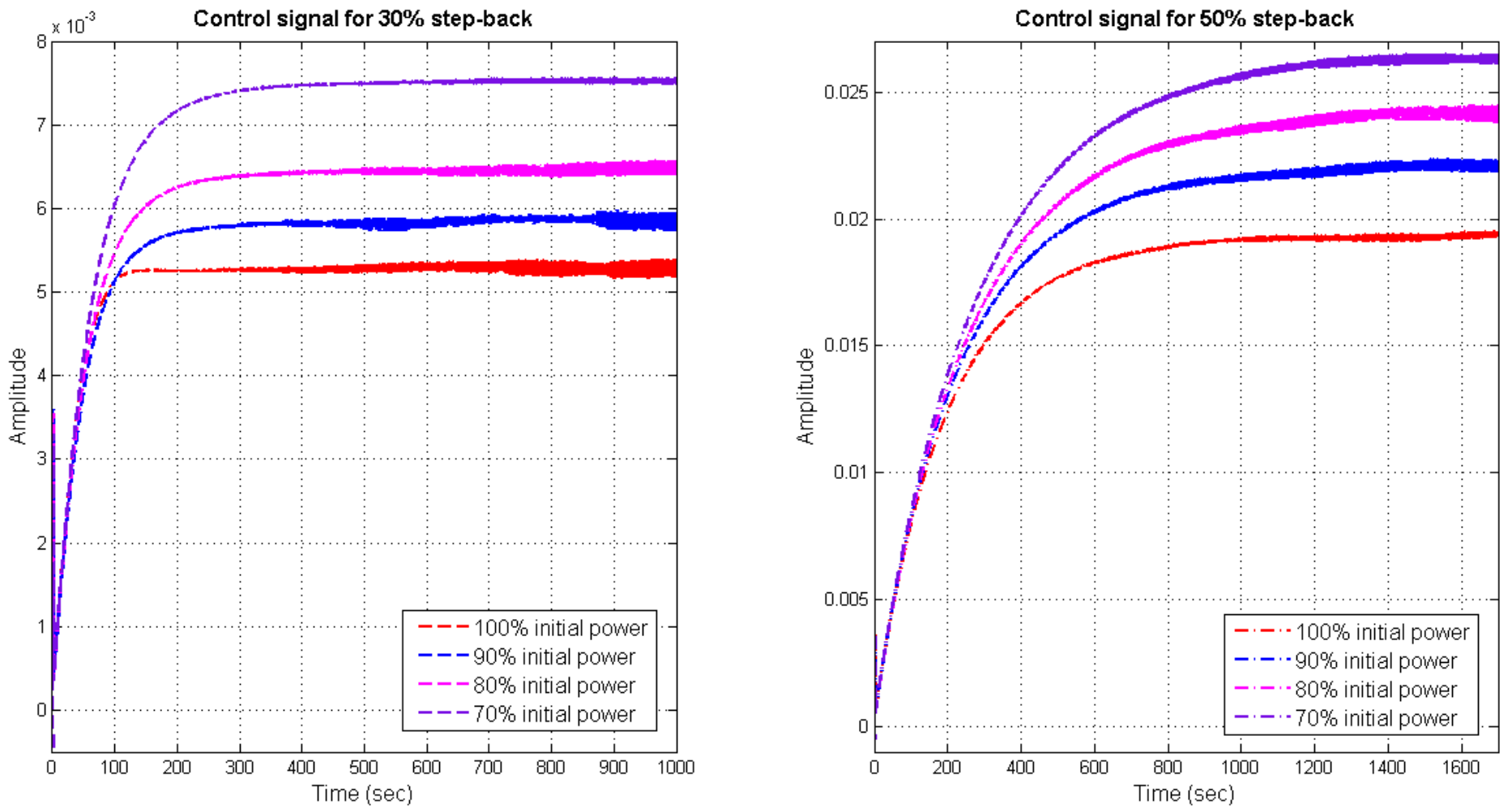

Fig. 17. Control signals for 30% and 50% step-back cases with unit step reference input.

In addition, the control signal or required variation in the control rod is shown in Fig. 17, for unit step reference input. It is observed from Fig. 17 that at lower initial powers, the required variation in control rod is higher due to decrease in the loop gain. Due to the same reason, larger control rod movement is needed for 50% step-back than that in the case of 30% step-back. The disturbance rejection performance of the designed controller is shown in Fig. 18. The disturbance rejection performance can be viewed as suppression of sudden reactivity inputs due to some other actions except rod movement and the capability of the controller to attenuate power oscillations due to such unwanted inputs. The simulations are reported with the models at 30% and 50% step-back subjected to a unit step disturbance input. It is observed that small local oscillations are present near the full power for 30% step-back due to high dc gain of those models.

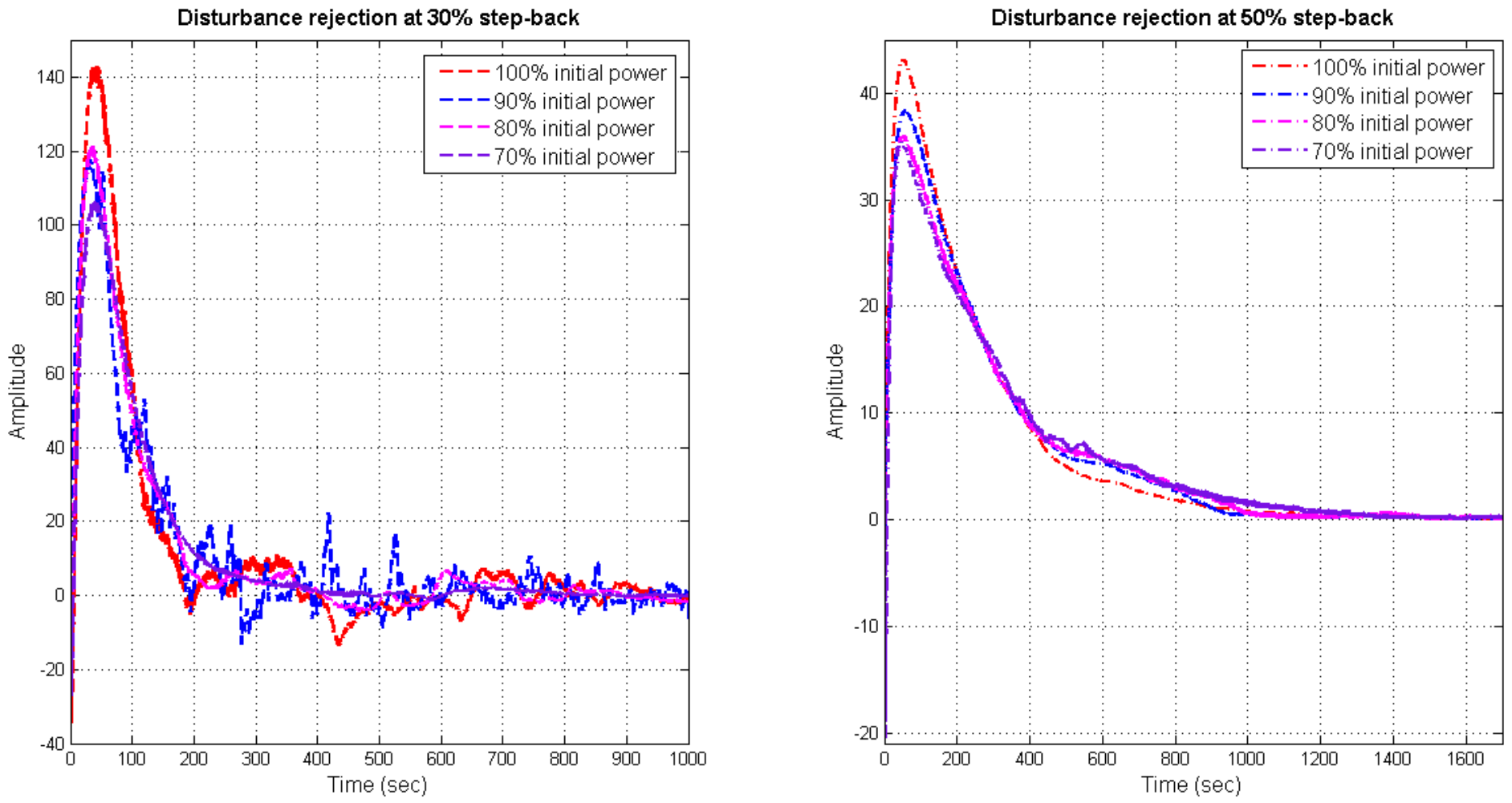

Fig. 18. Load disturbance rejection performance for 30% and 50% step-back model cases with unit step disturbance input.

### *4.3. Contributions of the design approach and few discussions*

The present approach looks the problem in a new way than that done in [8], [9]. PID control places the closed loop poles in the negative half of the s-plane i.e. the primary Riemann sheet. A FOPID controller as in [8] although being more robust than PID controller, may have few poles in higher Riemann sheets but still there will be few under-damped poles in the primary Riemann sheet. This may make the system to tend towards oscillations or instability when the gain of the plant increases excessively. Keeping in mind the danger under unpredictable catastrophic failures i.e. sudden and unusual increase in loop gain than the usual cases of operating the nonlinear system around different operating regimes (e.g. as in Das *et al.* [32]), the present approach focuses on placing the closed loop poles in the hyper-damped region in the complex $w$-plane for the fractional order system. Since the controller tuning algorithm drives all the closed loop poles in the secondary Riemann sheet, the chance of instability becomes almost insignificant. Even under tremendous gain increase, the closed loop pole has to cross the boundary between secondary and primary Riemann sheet, and thereafter cross the whole primary Riemann sheet before reaching marginal stability.

The enhanced safety issue with this approach of continuous order PID controller comes at the cost of slow reactor operation than that with the PID/FOPID controller in [8]. So, the focus of this work is to increase safety features with a new design philosophy and not to increase the performance of control. Similar studies on nuclear reactor power level control have been attempted in Das *et al.* [8] and the comparison of robustness between PID and FOPID controller has been shown. Since, in [8] performance comparison of PID/FOPID controller has already been done, we omitted similar comparisons from the present study. It is to be noted that beside the most essential performances like steady-state offset removal at all operating points, the main focus of the present work is to increase safety features which are obtained in the form of hyper-damping. In short the extra hyper-damped poles in close loop is making the system extra safe against very wide changes in system gain, which otherwise is not possible by conventional tuning of PID/FOPID controlled systems [8].

In addition, for power level adjustment of a nuclear reactor control rod movement is not the only means. In the point reactor kinetics the total reactivity may be changed in two ways viz. using the movement of control rods and changing the coolant temperature. The latter is known as thermal feedback or thermal-hydraulic effect on reactivity. Relevant detailed mathematical treatment and modeling have been reported in Das *et al.* [32]. In the present study we have only considered the rod movement for changing the reactivity levels, for the sake of simplicity. The rod movement is immediate action to correct the error, though additional shim controls, liquid zone control systems (LZCS) (for large PHWR) also exist. But in this paper the focus is towards primary device. Similar treatments may be put for the secondary fine control devices too like thermal feedback control, shim controls, LZCS etc.

Also, in order to design the continuous order controller to enjoy the safety features of hyper-damping, the number of zeros in different Riemann sheets ($N$) and the commensurate order ($q$) of the controller need to be fixed before tuning its gains using the proposed optimization based approach. The parameters of the controller (2) i.e. $N$ and q may be selected so as to match the $N$ and $q$ of the system under control which has a generic structure like (3). Firstly the maximum order of the controller ($Nq$) and commensurate order ($q$) are fixed by making them same as that of the system, so as to precisely move $N$ number of system poles using N number of controller zeros. In order to do that, an optimization based technique may be adopted to search for controller gains i.e. numerator coefficients while ensuring closed loop pole placement at desired locations. Use of conventional PID controllers may produce ample phase margin or over-damping at the cost of reduced performance but can never give hyper-damping. Thus to face nonlinearity and for added safety reasons hyper-damping with FO controllers is a better measure than wide phase margin or over-damping, since the latter may go to oscillation under violent increase in loop gain.

Every new design approach in order to improve reactor operation from performance or safety point of view is often questioned whether it's compliant with other constraints like maximum temperature decrease rate or not. From the simulations in Fig. 16, it is evident that the reactor is now being operated in quite slow rate compared to that reported in [8]. For faster reactor operation, temperature decrease rates are of big concern. So, with the proposed scheme, the heat removal mechanism is quite simpler to implement since the temperature increase or decrease is slower. We have not concentrated in the thermal corrections and by it final control which are usually done for large reactor for flux flattening purpose. Here our objective is to insert hyper-damped poles in different higher Riemann Sheets for the primary rod controller. This strategy has made the system slower compared to [8], but the thermal correction algorithms do follow the same. The thermal time constants are very large and thus not being considered for our investigation.

## 5. Conclusion

The paper reports a continuous order modeling approach for a nuclear reactor under varying step-back conditions in order to design a continuous order PID like controller. The data based reactor modeling is first attempted with four least square estimator variants to get discrete time transfer function models which have further been used to produce frequency domain data to build continuous order models with various levels of sampled order distribution. Frequency domain system identification technique is used to build the fractional order models with commensurate order 0.25 as a trade-off between complexity of the models and their accuracy. Optimization based pole assignment like approach has been adopted to design PID like continuous order controller in the $w$-plane having the same commensurate orders as the reactor models. The controller not only ensures dead-beat power level tracking at different operating conditions of the reactor but also ensures high reliability and safety at increased gain. The effectiveness of hyper-damped closed loop poles as design criterion ensures oscillation free power level tracking with enhanced stability as the root locus branches are far away from instability region due to increased loop gain caused by nonlinearity or possible mishandling by operator or in accidental condition.

Major findings of the present paper over the existing methodologies in continuous order system identification and controller design are as follows:

- Unnecessary refinement in the commensurate order for fractional order model building in order to achieve close approximation of continuous order model ($q \to 0$) may not be always beneficial, as the large system matrices become ill-conditioned and accuracy of the models decreases. So, an intuitive judgment is needed by looking at the commensurate order as well as the corresponding modeling accuracy to decide required refinement in sampled order distribution.
- The structure of continuous order controller has been chosen with several zeros in FO domain and a single integer order pole only. Introduction of such a controller would definitely increase the stability of the closed loop system as all the zeros attract the root locus branches and the single integrator works sufficiently well to eliminate the steady state off-set. In such cases, fractional order integrator with order less than unity can only be used if the designer can allow more sluggish time response.
- The order distribution curves in the contemporary literatures [10] shows monotonic increasing/decreasing nature or having some ideal and smooth distributions. The notion was to approximate the experimentally found order distribution with available curve fitting techniques to find out an equation representing the continuous order distribution. We found that for the reactor models the discrete order distributions are widely varying and also not in a regular manner. Therefore, concepts like finding dominant orders by just looking at the magnitude of the coefficients, finding the equation of the order distribution to get closed form expressions for continuous order transfer functions etc. cannot be applied under all circumstances. Still a sampled continuous order modeling based controller design in secondary Riemann sheet can be an effective way to design hyper-damped control systems for enhanced safety at high gains.
- The proposed hyper-damped controller design technique provides additional safety features against large gain variation in a faulty situation. The concept of hyper-damping, which can only be obtained using fractional order controllers, is especially useful in safety critical application like nuclear reactor power level maneuvering. Thus classical PID control loops are not capable of providing high robustness against large gain variation, which is the motivation of the present approach.

Future scope of research can be directed towards finding analytical closed form solution like in [10] for experimental data driven continuous order models with varying level of sampled order distribution and finding suitable control scheme to stabilize them.

**Acknowledgement:**
This work has been supported by Department of Science and Technology (DST), Govt. of India, under the PURSE programme.